\documentclass{article}

\usepackage{amsmath,amsfonts,amssymb,enumerate,amsthm}

\newtheorem{theorem}{Theorem}[section]
\newtheorem{definition}[theorem]{Definition}
\newtheorem{proposition}[theorem]{Proposition}
\newtheorem{corollary}[theorem]{Corollary}
\newtheorem{lemma}[theorem]{Lemma}
\newtheorem{remark}[theorem]{Remark}
\newtheorem{example}[theorem]{Example}

\newcommand{\bdfn}{\begin{definition}}
\newcommand{\edfn}{\end{definition}}
\newcommand{\bthm}{\begin{theorem} }
\newcommand{\ethm}{ \end{theorem}}
\newcommand{\bprop}{\begin{proposition}}
\newcommand{\eprop}{\end{proposition}}
\newcommand{\bcor}{\begin{corollary}}
\newcommand{\ecor}{\end{corollary}}
\newcommand{\blem}{\begin{lemma}}
\newcommand{\elem}{\end{lemma}}
\newcommand{\bfact}{\begin{remark} }
\newcommand{\efact}{\end{remark}}
\newcommand{\bex}{\begin{example} }
\newcommand{\eex}{ \end{example}}

\newcommand{\be}{\begin{enumerate}}
\newcommand{\ee}{\end{enumerate}}
\newcommand{\bce}{\begin{center}}
\newcommand{\ece}{\end{center}}
\newcommand{\bcev}{\vspace*{-0.2cm}\begin{center}}
\newcommand{\ecev}{\end{center}\vspace*{-0.2cm}}
\newcommand{\bi}{\begin{itemize}}
\newcommand{\ei}{\end{itemize}}
\newcommand{\bt}{\begin{tabular}}
\newcommand{\et}{\end{tabular}}
\newcommand{\vsp}[1]{\vspace*{#1cm}}

\newcommand{\ba}{\leftrightarrow}
\newcommand{\sra}{\rightarrow}
\newcommand{\ra}{\rightarrow}
\newcommand{\defi}{:=}
\newcommand{\bm}[1]{{\mathbf #1}}

\newcommand{\se}{\subseteq}
\newcommand{\vp}{\varphi}
\newcommand{\ti}[1]{\overline{#1}}

\newcommand{\real}{\mathbb R}
\newcommand{\fct}[1]{\widetilde{#1}}
\newcommand{\li}[1]{qspan(#1)}
\newcommand{\sq}{\varrho}
\newcommand{\rmvlog}{{\mathbb R}{\cal L}}
\DeclareMathOperator{\U}{U}

\DeclareMathOperator{\Gam}{\Gamma}
\DeclareMathOperator{\Gamr}{\Gamma_{\mathbb R}}

\title{\L ukasiewicz logic and Riesz spaces}

\author{Antonio Di Nola\\ Universita di Salerno,\\
 Via Ponte don Melillo  84084 Fisciano, Salerno, Italy\\ {\small adinola@unisa.it}\and 
Ioana Leu\c stean\\  University of Bucharest,\\Academiei 14, sector 1,  C.P. 010014, Bucharest, Romania\\ 
{\small ioana@fmi.unibuc.ro}}

\date{}
\begin{document}
\maketitle                  

\begin{abstract} We initiate a deep study of {\em  Riesz MV-algebras} which are  MV-algebras endowed with a scalar multiplication with scalars from $[0,1]$. Extending Mundici's equivalence between MV-algebras and $\ell$-groups, we prove that  Riesz MV-algebras are categorically equivalent with unit intervals in Riesz spaces with strong unit.   Moreover, 
the subclass of norm-complete Riesz MV-algebras is equivalent with the class of commutative unital C$^*$-algebras. 
The propositional calculus ${\mathbb R}{\cal L}$ that has Riesz MV-algebras as models is a conservative  extension of \L ukasiewicz $\infty$-valued propositional calculus and it is complete with respect to evaluations in the standard model $[0,1]$. We prove a normal form theorem for this logic, extending McNaughton theorem for \L  ukasiewicz logic.  
We define the notions of quasi-linear combination and quasi-linear span for formulas in ${\mathbb R}{\cal L}$ and we relate them with  the analogue of de Finetti's coherence criterion for ${\mathbb R}{\cal L}$.
\end{abstract}

\section{Introduction}

MV-algebras are the algebraic structures   of \L ukasiewicz $\infty$-valued logic.
The real unit interval $[0,1]$ equipped with the operations
\vspace*{-0.2cm}
\bce
$x^*=1-x$ and  $x\oplus y=\min(1,x+y)$ 
\ece
\vspace{-0.2cm}
for any $x, y\in [0,1]$, is the {standard} 
MV-algebra, i.e. an equation holds in any MV-algebra if and only if it holds in $[0,1]$.
In \cite{M1} Mundici  proved that
 MV-algebras are categorically equivalent with abelian lattice-ordered groups with strong unit. Consequently,
 for any MV-algebra there exists a lattice-ordered group  with strong unit $(G,u)$ such that $A\simeq [0,u]_G$, where
  \bcev
 $[0,u]_G=([0,u],\oplus,^*,0)$,\\
 $[0,u]=\{x\in G\mid 0\leq x\leq u\}$, \\
$x\oplus y=(x+y)\wedge u$ and $x^*=u-x$ for any $x$, $y\in [0,u]$.
\ece
If $(V,u)$ is a Riesz space (vector-lattice) \cite{LZ} with strong unit  then the
unit interval $[0,u]_V$ is closed to the scalar multiplication with scalars from $[0,1]$.  The structure 
\bcev $[0,u]_V=([0,u],\cdot,\oplus,^*,0)$, 
\ecev
where $([0,u],\oplus,^*,0)$ is the  MV-algebra defined as above  and $\cdot:[0,1]\times [0,u]_V\to [0,u]_V$  satisfies the axioms of the scalar product is the fundamental example in the theory of {\em Riesz MV-algebras},  initiated in \cite{DL} and further developed in the present paper.

The study of Riesz MV-algebras  is  related to the problem of
finding a complete axiomatization for the
variety generated by $([0,1],\cdot,\oplus,^*,0)$, where
$([0,1],\oplus,^*,0)$ is the standard MV-algebra and $\cdot$ is
the product of real numbers.
 The investigations led to the
definition of  {\em product MV-algebras} (PMV-algebras), which can be
represented as unit intervals  in lattice-ordered rings with strong
unit \cite{DiDv}.
A PMV-algebra is a structure $(P,\cdot)$, where $P$ is an MV-algebra and
$\cdot:P\times P\to P$ satisfies the equations of an internal product.
PMV-algebras are an equational class, but  the
 standard model $[0,1]$  generates only a quasi-variety which is
a proper subclass  of PMV-algebras \cite{Mo2}.  In this context, it was natural to replace the internal product with an external one:
a {\em Riesz MV-algebra} is a structure $(R,\cdot)$, where $R$ is an MV-algebra and $\cdot:[0,1]\times R\to R$.
Since we prove that the variety of Riesz MV-algebras is generated by $[0,1]$,  the  propositional calculus $\rmvlog$,  that has Riesz MV-algebras as models, is complete with respect to evaluations in $[0,1]$.

The study of Riesz MV-algebras was initiated in \cite{DL}.  In Section \ref{s:RMV} we give an equivalent but more suitable definition of these structures  and we prove some of their fundamental  properties.

The categorical equivalence between  Riesz MV-alge-bras and Riesz spaces with strong unit is proved in Section \ref{comp}.
 As a consequence,  the standard Riesz MV-algebra $[0,1]$ generates the variety of Riesz MV-algebras. 

In Section \ref{sec-norm}, the categorical equivalence is specialized to the class of {\em norm-complete Riesz MV-algebras}, which is  dually equivalent with the category of compact Hausdorff spaces. Using the Gelfand-Naimark duality, this leads us to a connection with the theory of commutative unital C$^*$-algebras.

Section \ref{propc}  presents the propositional calculus $\rmvlog$ which simplifies the one introduced in \cite{DL}.
In Section \ref{normal} we prove a normal form theorem for formulas of  $\rmvlog$. Since $\rmvlog$
is a conservative extension of \L ukasiewicz logic $\cal L$, this theorem is a generalization of
McNaughton theorem \cite{McN}. Our result asserts that
{\em for any continuous piecewise linear  function $f:[0,1]^n\to [0,1]$  there exists a formula $\varphi$ of
$\rmvlog$ with $n$ variables such that $f$ is the term function associated to $\varphi$}.

In Section \ref{linear} we  initiate the theory of{\em quasi-linear combinations} of formulas in  $\rmvlog$.
If $f_i:[0,1]^n\to{\mathbb R}$ are continuous piecewise linear
 functions
and  $c_i$ are real numbers for any $i\in \{1,\ldots,k\}$, then the normal form theorem guarantees the existence of a formula $\Phi$ of $\rmvlog$, whose term function is equal to $((\sum_{i=1}^kc_if_i)\vee 0)\wedge 1$ and, in this case, we say that $\Phi$ is {\em quasi-linear combination} of  $f_1$, $\cdots$, $f_k$. We prove de Finetti's coherence criterion for $\rmvlog$  and we 
 provide an equivalent characterization by the fact that a  quasi-linear span contains only
invalid formulas. 

Some of results contained in this paper may overlap with the proceeding paper \cite{DL}. Other results are proved in  a more general setting in \cite{DFL,FL1,Lbook}.  For the sake of completeness we sketched the proofs that we consider important for the present development.   
 
\section{Preliminaries  on MV-algebras}

An {\em MV-algebra}  is a structure $(A,\oplus,^{*},0)$, where
$(A,\oplus, 0)$ is an abelian monoid and the following identities hold
for all $x,y\in A$:
\begin{enumerate}[(MV1)]
\item $(x^{*})^{*}=x$,
\item $0^{*}\oplus x=0^{*}$, 
 \item $(x^{*}\oplus y)^{*}\oplus y=(y^{*}\oplus x)^{*}\oplus x$.
\end{enumerate}
 We refer to \cite{CDOM} for all the unexplained notions concerning MV-algebras and to \cite{DM} for advanced topics. On any MV-algebra $A$ the following operations are defined for any $x,y\in A$:
\bce
$1=0^*$, $x\odot y=(x^*\oplus y^*)^*$, $x\sra y=x^*\oplus y$

$0x=0$, $mx=(m-1)x\oplus x$ for any $m\geq 1$.
\ece
We assume in the sequel that the operation $\odot$ is more binding then $\oplus$.

\bfact\label{ordmv}
Any  MV-algebra $A$ is a bounded distributive lattice, with  the partial order defined by
\bcev
$x\leq y$ if and only if $x\odot y^*=0$
\ecev
and the lattice operations defined by 
\bcev
$x\vee y=x\oplus y\odot x^*$ and $x\wedge y=x\odot(x^*\oplus y)$
\ecev
for any $x$, $y\in A$.
\efact

Any MV-algebra $A$  has an internal  distance:
\bcev
$d(x,y)=(x\odot y^*)\oplus (x^*\odot y)$ for any $x$, $y\in A$.
\ecev

\blem\label{dist}
\cite[Proposition 1.2.5]{CDOM} In any MV-algebra $A$, the following properties hold for any $x$, $y$, $z\in A$:\\
(a) $d(x,y)=d(y,x)$,\\
(b) $d(x,y)=0$ iff $x=y$,\\
(c) $d(x,z)\leq d(x,y)\oplus d(y,z)$.
\elem

If $(A,\oplus,^*,0)$ is an MV-algebra then an {\em ideal} is a nonempty subset $I\se A$ such that for any $x$, $y\in A$ the following conditions are satisfied: 
\begin{enumerate}
\item[(i1)] $x\in I$ and  $y\leq x$ imply $y\in I$,
\item[(i2)] $x$ and  $y\in I$ imply $x\oplus y\in I$.
\end{enumerate}
An ideal $I$ of $A$ uniquely defines a congruence $\sim_I$ by
\bcev
$x\sim_I y$ iff $x\odot y^*\in I$ and $y\odot x^*\in I$. 
\ecev 
We denote by $A/I$ the quotient MV-algebra and we refer to \cite{CDOM} for more details.

We recall that an {\em $\ell$-group}  is a structure $(G,+,0,\leq)$ such that $(G,+,0)$ is a group, $(G,\leq)$ is a lattice and  any group translation is isotone \cite{BKW}. In the following the $\ell$-groups are abelian.
For an $\ell$-group $G$ we denote $G_+=\{x\in G\mid x\geq 0\}$.
An element $u\in G$ is  a {\em strong unit} if $u\geq 0$ and for any $x\in G$ there is a natural number $n$ such that $x\leq nu$. An $\ell$-group  is {\em unital} if it posses a strong unit. 
If $(G,u)$ is a unital $\ell$-group,  we define  $[0,u]=\{x\in G\mid 0\leq x\leq u\}$ and
\bce
$x\oplus y=(x+y)\wedge u$,  $ x^*=u-x$   for any $x, y\in [0,u]$.
\end{center}
Then $[0,u]_G=([0,u],\oplus,\neg,0)$ is an MV-algebra \cite[Proposition 2.1.2]{CDOM}.

\begin{lemma}\label{l:goodsequence}\cite[Lemma 7.1.3]{CDOM}
Let  $(G,u)$ be a  unital $\ell$-group,  $x\geq 0$ in $G$ and $n\geq 1$ a natural number such that $x\leq nu$. Then 
$x=x_1+\cdots + x_n$,
where
 \bce $x_i=((x-(i-1)u)\vee 0)\wedge u\in [0,u]$\ece for any $i\in\{1,\ldots,n\}$.
\end{lemma}

We denote by ${\mathbb M}{\mathbb V}$  the category of MV-algebras and by ${\mathbb A}{\mathbb G}_u$  the category of unital abelian lattice-ordered groups with unit-preserving morphisms.
  In \cite{Meq} the functor $\Gam\colon {\mathbb A}{\mathbb G}_u\to{\mathbb M}{\mathbb V}$   is defined as follows:\\
$\Gam(G,u)=[0,u]_G$ for any unital $\ell$-group $(G,u)$, \\
$\Gam(f)=f|_{[0,u]}$ for any  morphsim $f:(G,u)\to (G^{\prime},u^\prime)$ from ${\mathbb A}{\mathbb G}_u$.
\begin{theorem}\label{t:mundici}\cite[Corollary 7.1.8]{CDOM}
 The functor $\Gam$ yields an equivalence between  ${\mathbb A}{\mathbb G}_u$ and ${\mathbb M}{\mathbb V}$. 
\end{theorem}

\bdfn
If $A$ and $B$ are  MV-algebras then a function $\omega:A\to B$ is called {\em additive} if 
\bce
$x\odot y=0$ implies $\omega(x)\odot\omega(y)=0$ and $\omega(x\oplus y)=\omega(x)\oplus\omega(y)$.
\ece
\edfn

Additivity was firstly studied in the context of states  defined on MV-algebras \cite{M1}. The theory of states generalizes the  boolean probability theory and reflects the theory of states defined on $\ell$-groups.

\bdfn\cite{M1}
If $A$ is an MV-algebra then a  function  $s:A\to [0,1]$ is a {\em state}  if the following properties are satisfied for
any  $x$, $y\in A$:
\begin{enumerate}[(s1)]
\item if  $x\odot y=0$ then 
$s(x\oplus y)=s(x)+s(y)$,
\item  $s(1)=1$.
\end{enumerate}
\edfn

The following results are proved in \cite{Lbook}, but we sketch the proofs for the sake of completeness. We also note that particular instances of these results are proved in \cite{FL1}. Proposition \ref{extop} is proved for states in \cite{M1}.

\bprop\label{extop}
Assume $(G,u)$ and $(H,v)$ are unital $\ell$-groups, $A=\Gam(G,u)$ and $B=\Gam(H,v)$. Then 
for any additive function $\omega:A\to B$ there exists a unique group morphism $\ti{\omega}:G\to H$ such that 
$\ti{\omega}(x)=\omega(x)$ for any $x\in [0,u]$.
\eprop
\begin{proof}
If $x\in G$ and $x\geq 0$ then there are $x_{1}$, $\ldots$, $x_{m}\in [0,u]$ such that 
$x=x_{1}+\cdots+x_{m}$. Then we define
\vsp{-0.2}
\bce
$\ti{\omega}(x)\defi \omega(x_{1})+\cdots +\omega(x_{m})$.
\ece
\vsp{-0.2}
The fact that  $\ti{\omega}(x)$ is well defined follows by Riesz decomposition property in $\ell$-groups \cite[1.2.16]{BKW}.
Hence $\ti{\omega}(x)$ is well defined for $x\in G_{+}$ and 
$\ti{\omega}(x+y)=\ti{\omega}(x)+\ti{\omega}(y)$ 
for any $x$, $y\in G_{+}$.
By \cite[1.1.7]{BKW} it follows  that $\ti{\omega}$ can be uniquely
extended to a group homomorphism defined on $G$.    
\end{proof}

\blem\label{add}
If $A$ and $B$  are  MV-algebras and $\omega:A\sra B$ is a function, then the following are equivalent:\\
(a) $\omega$ is additive,\\
(b) the following properties hold  for any $x$, $y\in A$:

(b1) $x\leq y$ implies $\omega(x)\leq\omega(y)$,

 (b2) $\omega(x\odot(x\wedge y)^{*})=\omega(x)\odot\omega(x\wedge y)^{*}$.
\elem
\begin{proof} 
(a)$\Rightarrow$(b) 
If $x\leq y$ then $y=x\vee y=x\oplus y\odot x^*$, so  $\omega(y)=\omega(x)\oplus \omega(y\odot x^*)$ and 
$\omega(x)\leq\omega(y)$. Hence, $\omega$ is isotone.
 We remark that $\omega(x\wedge y)\leq \omega(x)$, so
\bcev
$\omega(x)\odot\omega(x\wedge y)^{*}\oplus \omega(x\wedge y) = \omega(x)\vee\omega(x\wedge y)
 =\omega(x) =\omega(x\vee(x\wedge y))
 =\omega(x\odot(x\wedge y)^{*})\oplus \omega(x\wedge y).$
\ecev
It follows that
\bcev
$\begin{array}{l}
\omega(x)\odot\omega(x\wedge y)^{*}=\\
(\omega(x)\odot\omega(x\wedge y)^{*})\wedge \omega(x\wedge y)^*=\\
(\omega(x)\odot\omega(x\wedge y)^{*}\oplus \omega(x\wedge y))\odot \omega(x\wedge y)^*=\\
(\omega(x\odot(x\wedge y)^{*})\oplus \omega(x\wedge y))\odot \omega(x\wedge y)^*=\\
\omega(x\odot(x\wedge y)^{*})\wedge \omega(x\wedge y)^*=\\
\omega(x\odot(x\wedge y)^{*})
\end{array}$
\ecev
 
\noindent (b)$\Rightarrow$(a) We remark that for $x=1$ in (b2) we get 
$\omega(y^{*})=\omega(1)\odot\omega(y)^{*}$, so \mbox{$\omega(y^{*})\leq\omega(y)^{*}$}.
Assume  $x\odot y=0$, so $x\leq y^{*}$. Using (b1), we get $\omega(x)\leq\omega(y^{*})\leq\omega(y)^{*}$, so $\omega(x)\odot \omega(y)=0$.
In this case, using (b2) we get
\vsp{-0.2}
\bce
$\omega(x)=\omega(x\wedge y^{*})=\omega((x\oplus y)\odot y^{*})=\omega(x\oplus y)\odot\omega(y)^{*}$.
\ece
\vsp{-0.2}
It follows that:
\vsp{-0.2}
\bce
 $\omega(x)\oplus \omega(y)=\omega(x\oplus y)\odot\omega(y)^{*}\oplus \omega(y)=\omega(x\oplus y)\vee\omega(y)$.
\ece
\vsp{-0.2}
Using (b1), $\omega(y)\leq\omega(x\oplus y)$, and we get $\omega(x\oplus y)=\omega(x)\oplus \omega(y)$.  
\end{proof}

\section{Riesz MV-algebras}\label{s:RMV}

 Riesz MV-algebras are introduced in \cite{DL}.  Below we give a 
simpler and more suitable definition, which provides directly an equational characterization. The equivalence between  this definition and the one from \cite{DL} is proved in \mbox{Theorem \ref{ech}.}

\bdfn\label{rmv.1}
A {\em Riesz MV-algebra} is a structure
\bce 
$(R, \cdot,\oplus,^*,0)$,
\ece
where $(R,\oplus,^*,0)$ is an MV-algebra and the operation $\cdot:[0,1]\times R\ra R$ satisfies the following identities  for any $r$, $q\in [0,1]$ and $x$, $ y\in R$:
\begin{enumerate}[(RMV1)]
\item\label{r1}  $r\cdot (x\odot y^{*})=(r\cdot  x)\odot(r\cdot y)^{*}$,
\item\label{r2} $(r\odot q^{*})\cdot x=(r\cdot x)\odot(q\cdot x)^{*}$,
\item\label{r3} $r\cdot (q\cdot  x)=(rq)\cdot x$,
\item\label{r4}  $1\cdot x=x$.
\end{enumerate}

\medskip

In the following we  write
$rx$ instead of $r\cdot x$ for $r\in [0,1]$ and $x\in R$. Note that $rq$ is the real product for any $r$, $q\in [0,1]$.
\edfn

\bex
If $X$ is a compact Hausdorff space then 
\bcev $C(X)_u=\{f:X\to [0,1]\mid  f\,\, \mbox{continuous}\}$
\ecev 
 is a Riesz MV-algebra, with all the operations defined componentwise. This example will be further investigated in Section \ref{sec-norm}
\eex

\bex
If $G$ is an abelian $\ell$-group, then $R=\Gam({\mathbb R}\times_{lex} G, (1,0))$  is a Riesz MV-algebra, where 
${\mathbb R}\times_{lex} G$ is the lexicographic product of $\ell$-groups and the scalar multiplication is defined by 
$r(q,x)=(rq, x)$ for any $r\in [0,1]$ and $(q,x)\in R$.

\eex

\blem\label{rmv-pr}
In any Riesz MV-algebra $R$  the following properties hold for any $r$, $q\in [0,1]$ and $x$, $y\in R$:\\
(a) $0x=0$, $r0=0$,\\
(b) $x\leq y$ implies $rx\leq ry$,\\
(c) $r\leq q$ implies $rx\leq qx$,\\
(d) $rx\leq x$.

\elem
\begin{proof}
(a)  follows by (RMV1) and (RMV2) for $x=y$ and, respectively, $r=q$. \\
(b), (c) follow by Remark \ref{ordmv}.\\
(d)  follows by (c) and (RMV4).   

\end{proof}

\bprop\label{embed}
 The function
$\iota :[0,1]\to R$ defined by $\iota(r)=r1$ for any $r\in [0,1]$
 is an embedding. Consequently, 
any Riesz MV-algebra $R$  contains a subalgebra isomorphic with $[0,1]$.
\eprop
\begin{proof}
By Lemma \ref{rmv-pr} we get $\iota(0)=0$. If $r$, $q\in [0,1]$ then 
\bce
$\iota(r^*)=r^* 1=(1\cdot 1)\odot (r1)^*=(r1)^*$,

$\iota(r\odot q)=\iota(r\odot q^{**})= (r\odot q^{**})1=(r1)\odot (q^* 1)^*=(r1)\odot (q1)^{**}=(r1)\odot (q1)$.  
\ece   

\end{proof}

A {\em Riesz space} ({\em vector lattice}) \cite{LZ} is a structure 
\bcev $(V,\cdot, +,0,\leq)$
\ecev
 such that $(V,+,0,\leq)$ is an abelian $\ell$-group, $(V,\cdot,+,0)$ is a real vector space and, in addition,\\
(RS) $x\leq y$ implies $r \cdot x\leq r\cdot y$ ,\\
for any $x$, $y\in V$ and $r\in{\mathbb R}$, $r\geq 0$.\\
A Riesz space is {\em unital} if the underlaying  $\ell$-group  is {unital}.

\blem\label{ex1}  If $(V,u)$ is a unital Riesz space, then
\bcev
$[0,u]_V=([0,u],\cdot,\oplus, ^*,0)$ 
\ecev
is a Riesz MV-algebra, where $rx$ is the scalar multiplication of $V$ for any $r\in [0,1]$ and 
$x\in [0,u]$.  
\elem
\begin{proof} Assume $r, q\in [0,1]$ and $x, y\in [0,u]$.\\
(RMV1) $r(x\odot y^*)=r((x-y)\vee 0)=(r x
-r y)\vee 0=(r x)\odot (r y)^*$.\\
(RMV2) If $r\leq q$ then $r x \leq q x$, so
$(r\odot q^*)x=$\\ $((r-q)\vee 0)x=0=(rx-qx)\vee 0=(r x)\odot ( q x)^*$.

If $r> q$ then
$((r-q)\vee 0)x=(r-q)x=rx-q x =(r x)-(q x)\vee 0=(r x)\odot( q x)^*$.

We note that (RMV\ref{r3}) and (RMV\ref{r4}) hold in $V$, therefore they hold in $[0,u]$.  
\end{proof}

\bfact
If  $(R, \cdot,\oplus,^*,0)$  is a Riesz MV-algebra then  we denote its MV-algebra reduct  by  $\U (R)=(R,\oplus,^*,0)$.
Assume $I$ is an ideal of $\U (R)$. By Lemma \ref{rmv-pr} (d) we infer that 
$rx\in I$ whenever $r\in [0,1]$ and $x\in I$. It follows, by (RMV1),  that 
$ rx\sim_I ry $ whenever $r\in [0,1]$ and $x\sim_I y$. 
 As consequence, the quotient $R/I$ has a canonical structure of Riesz MV-algebra.
\efact

\bfact\label{fact-id}
A Riesz MV-algebra $R$ has the same theory of ideals (congruences) as its reduct $\U (R)$. 
If $R$ is a Riesz MV-algebra and $P\se R$ an ideal then it is straightforward that the following hold:\\
(a)  $P$ is prime iff $R/P$ is linearly ordered,\\
(b) $P$ is maximal iff $R/P\simeq [0,1]$. \\
Note that (b) holds since, for any maximal ideal $P$,  the quotient $R/P$ is an MV-subalgebra of $[0,1]$. But the only subalgebra of $[0,1]$ which is a Riesz MV-algebra is $[0,1]$ by Proposition \ref{embed}, so $R/P\simeq [0,1]$. 
\efact

\blem\label{help1} 
If $R$ is a Riesz MV-algebra, $I\se R$ an ideal and $x\in R$ such that $rx\in I$ for some 
$r\in (0,1]$ then $x\in I$.
\elem
\begin{proof} Let $r\in (0,1]$ such that $rx\in I$ and let $m$ be the integer part of $\frac{1}{r}$. Hence $\frac{1}{m+1}x\leq rx$, so $\frac{1}{m+1}x\in I$. Since $x=(m+1)(\frac{1}{m+1}x)$ we get $x\in I$. 
 
\end{proof}

\bcor
Any simple Riesz MV-algebra is isomorphic with $[0,1]$.  Any semisimple Riesz MV-algebra is a subdirect product of copies of $[0,1]$.
\ecor

In the sequel we investigate the  morphisms of Riesz MV-algebras.

\bcor\label{morf}
If $R_1$ and $R_2$ are Riesz MV-algebras and $f:U(R_1)\to U(R_2)$ is a morphism of MV-algebras then
\bce
$f(rx)=rf(x)$ for any $r\in [0,1]$ and $x\in R_1$. 
\ece
\ecor
\begin{proof}
 Assume $J$ is an ideal in $R_2$. Since $f$ is an morphism of MV-algebras, it follows that $f^{-1}(J)$ is an ideal in $R_1$.
If  $x\in R_1$ and $r\in [0,1]$  we have

\begin{center}
$r f(x)\in J\Rightarrow  f(x)\in J\Rightarrow x\in f^{-1}(J)
\Rightarrow$ 

$ r x \in f^{-1}(J)\Rightarrow f(r x)\in J$,

$f(r x)\in J\Rightarrow r x \in f^{-1}(J) 
\Rightarrow x \in f^{-1}(J) \Rightarrow$

$  f(x) \in J \Rightarrow r f(x) \in J$.
\end{center}

Note that we used Lemma \ref{help1} twice. We proved that, for any ideal $J$ of $R_2$
\bce
$r f(x)\in J \Leftrightarrow f(r x)\in J.$
\ece

Therefore $rf(x)\odot f(rx)^*\in J $ and $f(rx)\odot(rf(x))^*$ for any ideal $J$ of $R_2$. This means that
$rf(x)\odot f(rx)^*=f(rx)\odot(rf(x))^*=0$, so $f(rx)=rf(x)$.   
\end{proof}

\bfact\label{fact-mor}
The above result asserts that a morphism of Riesz MV-algebra is simply a morphism between the corresponding MV-algebra reducts. 
\efact

The following result is similar with  Chang's representation theorem for MV-algebras \cite[Lemma 3]{chang}.

\bcor\label{rep}
Any Riesz MV-algebra is  a subdirect product of linearly ordered Riesz MV-algebras. 
\ecor
\begin{proof}
If $R$ is a Riesz MV-algebra then, by Remark \ref{fact-id}, $\bigcap\{P\mid P\mbox{ prime ideal of } R\}=\{0\}$ and   $R/P$ is linearly ordered for any prime ideal $P$. As consequence, $R$ is a subdirect product 
of the family 
\bce $\{R/P\mid P\mbox{ prime ideal of } R\}$.  
 \ece
\end{proof}

In order to prove that Riesz MV-algebras introduced in Definition \ref{rmv.1} coincide with the ones defined in \cite{DL}, we recall some results from \cite{FL1}.

\bfact\cite{FL1}\label{fl}
If  $\Omega$ is a set of  unary operation symbols,  then an {\em MV-algebra with  $\Omega$-operators} is a structure $(A,\Omega_A)$ where  $A$ is an MV-algebra and for any $\omega\in \Omega$ the operation $\omega_A:A\to A$ is additive.  
An additive  function 
$\omega : A\to A$ is an {\em $f$-operator} if
\bce  $x\wedge y=0$ implies $\omega(x)\wedge y=0$ for any $x$, $y\in A$.
\ece

If $(A,\Omega_A)$ is an MV-algebra with $\Omega$-operators such that $\omega_A$ is an $f$-operator for any 
$\omega\in \Omega$, then  $(A,\Omega_A)$ is a subdirect product of linearly ordered MV-algebras with $\Omega$-operators 
\cite[Corollary 5.6. ]{FL1}. 
\efact

\bfact\label{fl1}
Assume that  $(R,\oplus,^*,0)$ is an MV-algebra and let  $\cdot:[0,1]\times R\to R$ such that 
{\rm (RMV2)}, {\rm (RMV3)}, {\rm (RMV4)} hold and the  function
\bcev
$\omega_r:R\to R$, $\omega_r(x)=r\cdot x$     
 \ecev
is additive for any $r\in [0,1]$. 
By {\rm (RMV2)} and  {\rm (RMV4)} we get  $\omega_r(x)\leq x$ for 
any $r\in [0,1]$ and $x\in R$, so   $\omega_r$ is an $f$-operator for any $r\in [0,1]$. If $\Omega=\{\omega_r\mid r\in [0,1]\}$ then,   by Remark \ref{fl},
 $(R,\Omega)$ is an MV-algebra with $\Omega$-operators  that can be represented  as 
 subdirect product of linearly ordered MV-algebras with $\Omega$-operators. 
\efact

\blem\label{help}
Assume that   $(R,\oplus,^*,0)$ is an MV-algebra. If $\cdot:[0,1]\times R\to R$ then the following are equivalent:\\
(RMV2) $(r\odot q^{*})\cdot x=(r\cdot x)\odot(q\cdot x)^{*}$ 

for any $r$, $q\in [0,1]$ and $x\in R$,\\
{\rm (RMV2$^\prime$)}  $r\odot q=0$ then $(r\cdot x)\odot (q\cdot x)=0$ and 
\bcev
$(r\oplus q)\cdot x=(r\cdot x)\oplus(q\cdot x)$ 
\ecev

for any $r$, $q\in [0,1]$ and $x\in R$.
\elem
\begin{proof}
 For $x\in R$ define $\omega_x: [0,1]\to R$ by \mbox{$\omega_x(r)=rx$} for any $r\in [0,1]$.
If $\omega_x$ satisfies (RMV2) then  the condition (b) from Lemma \ref{add} is satisfied, so $\omega_x$ satisfies also (RMV2$^\prime$).
Conversely, if $\omega_x$ satisfies (RMV2$^\prime$)  then, by Lemma \ref{add}, we also get $\omega_x(0)=0$.   Assume $r$, $q\in [0,1]$ such that  $r\leq q$.  Hence  $r\odot q^*=0$ and $rx\leq qx$, so  $ (r\odot q^*)x=0x=0=(rx)\odot (qx)^*$. If 
$r$, $q\in [0,1]$ such that $r>q$ then  (RMV2) coincides with the equation (b2) from Lemma \ref{add}.  
\end{proof}

\bthm\label{ech}
Assume that  $(R,\oplus,^*,0)$ is  an MV-algebra and  $\cdot:[0,1]\times R\to R$. Then   $(R,\cdot, \oplus,^*,0)$
 is a Riesz MV-algebra if and only if the following properties  are satisfied for any $x$, $y\in R$ and $r$, $q\in [0,1]$:\\
{\rm (RMV1$^\prime$)} if $x\odot y=0$ then $(r\cdot x)\odot(r\cdot y)=0$ and 
\bcev
$r\cdot (x\oplus y)=(r\cdot x)\oplus(r\cdot y)$,
\ecev
{\rm (RMV2$^\prime$)} if  $r\odot q=0$ then $(r\cdot x)\odot (q\cdot x)=0$ and 
\bcev
$(r\oplus q)\cdot x=(r\cdot x)\oplus(q\cdot x)$,
\ecev
{\rm (RMV3)} $r\cdot (q\cdot  x)=(rq)\cdot x$,\\
{\rm (RMV4)} $1\cdot x=x$.
\ethm
\begin{proof}
 By Lemma \ref{help},
if  $(R,\oplus,^*,0)$ is an MV-algebra and  $\cdot:[0,1]\times R\to R$, then   the algebra $(R,\cdot, \oplus,^*,0)$  satisfies (RMV2), (RMV3) and (RMV4) if and only if it satisfies  (RMV2$^\prime$), (RMV3) and (RMV4). 
Assume now that  $(R,\cdot, \oplus,^*,0)$  satisfies (RMV2), (RMV3) and (RMV4). We have to prove that (RMV1) is satisfied if and only if (RMV1$^\prime$) is satisfied. By Corollary \ref{rep} and  Remark \ref{fl1},  it suffices to prove the equivalence  for linearly ordered structures.  In this case, by Lemma \ref{add}, the equivalence of (RMV1) and (RMV1$^\prime$) is straightforward.    
 
\end{proof}

 Note that in \cite{DL} a Riesz MV-algebra is defined by  (RMV1$^\prime$), (RMV2$^\prime$), (RMV3) and (RMV4), so we proved that 
 Definition \ref{rmv.1} is equivalent with the initial one.

\section{Riesz MV-algebras and Riesz spaces. The completeness theorem.}\label{comp}

By Theorem \ref{ech}, Riesz MV-algebras are exactly the MV-modules \cite{DFL} over $[0,1]$.  Hence some basic properties follow from the general theory of MV-modules developed in \cite{DFL,Lbook}.  One of the most important results is the categorical equivalence between Riesz MV-algebras and unital Riesz spaces.  For the sake of completeness, we sketch the proof of this result.

\bprop\label{rep}
For any Riesz MV-algebra $R$ there is a unital Riesz space $(V,u)$ such that $R\simeq [0,u]_V$.
\eprop
\begin{proof}

By Theorem \ref{t:mundici},  there exists a unital $\ell$-group  $(V,u)$  such that $R$ and $[0,u]_V$ are isomorphic MV-algebras.
For any $\lambda\in {\mathbb R}$ and $x\in V$ we have to define  the scalar multiplication $\lambda x$.
We can safely assume that $R=[0,u]\subseteq V$.

If $r\in [0,1]$  then $x\mapsto r x$ is an additive function from $[0,u]_V$ to $[0,u]_V$ so, by Proposition 
\ref{extop}, it can be uniquely extended to a group morphism $\omega_r:V\to V$. Hence we define
$r x=\omega_r(x)$ for any $x\in V$. We note that $x\geq 0$ implies $r x\geq 0$.

If $q\in [0,1]$ then $\omega_{rq}=\omega_r\circ\omega_q$ since they coincide on the positive cone, so $r(q v)=(rq)v$. 

Note that $v=(v\vee 0)-((-v)\vee 0)$, so
\bcev
 $1v=1(v\vee 0)-1((-v)\vee 0)=(v\vee 0)-((-v)\vee 0)=v$.
\ecev
If  $\lambda\geq 0$ and $v\in V$, then 
there are  $r_{1}$, $\ldots$, $r_{m}\in [0,1]$ such that
$\lambda=r_{1}+\cdots+r_{m}$. Then we define
\bcev
$\lambda v = r_1 v+\cdots+r_m v.$
\ecev
One can prove  that $\lambda v$ is well-defined using the Riesz decomposition property \cite[1.2.16]{BKW}.   
If $\mu\geq 0$, then $\mu=q_{1}+\cdots+q_{n}$ for some $q_{1}$, $\ldots$, $q_{n}\in [0,1]$ and 
\begin{center}
$\lambda(\mu v)=\lambda(\sum_{j=1}^{n}q_j v)=\sum_{i=1}^{m}r_i\left(\sum_{j=1}^{n}q_j v\right)=\sum_{i=1}^{m}\sum_{j=1}^{n}r_i(q_j v)=
\sum_{i=1}^{m}\sum_{j=1}^{n}(r_iq_ j)v=
\left(\sum_{i=1}^{m}\sum_{j=1}^{n}(r_iq_ j)\right)v=(\lambda\mu) v$.
\end{center}

If $\lambda\leq 0$ in ${\mathbb R}$ then we set $\lambda v=-(|\lambda|v))$, where $|\lambda|$ is the module of $\lambda$ in ${\mathbb R}$. 
It is straightforward that $\lambda(\mu v)=(\lambda\mu) v$ for another $\mu\in {\mathbb R}$.

We know that $(V,u)$ is a unital $\ell$-group and we defined the scalar product $\lambda v$ for any $\lambda\in R$ and $v\in V$ such that $\lambda v\geq 0$  whenever $\lambda\geq 0$ and $v\geq 0$. 
 Therefore, $(V,u)$ is a unital vector lattice.
 
\end{proof}

We denote by $ {\mathbb R}{\mathbb M}{\mathbb V}$   the category of Riesz MV-algebras and by ${\mathbb R}{\mathbb S}_u$  the category of unital Riesz spaces with unit-preserving morphisms.

Following this construction we get a functor \bce$\Gamr\colon {\mathbb R}{\mathbb S}_u\to{\mathbb R}{\mathbb M}{\mathbb V}$\ece defined as follows:\\
$\Gamr(V,u)=[0,u]_V$ for any unital Riesz space $(V,u)$, \\
$\Gamr(f)=f|_{[0,u]}$ for any  morphism  $f:(V,u)\to (V^{\prime},u^\prime)$ from  ${\mathbb R}{\mathbb S}_u$.

\bthm\cite{DL,DFL}\label{cat}
The functor $\Gamr$ yields an equivalence between  ${\mathbb R}{\mathbb S}_u$ and ${\mathbb R}{\mathbb M}{\mathbb V}$. 
\ethm
\begin{proof}
It follows from Theorem \ref{t:mundici}, Corollary \ref{morf} and Proposition \ref{rep}.  
\end{proof}

It is straightforward that the following diagram is commutative, where $\U$ are forgetful functors:
\bce
\begin{tabular}{ccc}
$ {\mathbb R}{\mathbb S}_u$& $\stackrel{\Gamr}{\longrightarrow}$&${\mathbb R}{\mathbb M}{\mathbb V}$\\
$\U\downarrow $&  &$\downarrow\U $\\
$ {\mathbb A}{\mathbb B}_u$& $\stackrel{\Gam}{\longrightarrow}$&${\mathbb M}{\mathbb V}$
\end{tabular}
\ece
\medskip

The standard Riesz MV-algebra is $([0,1],\cdot,\oplus,^*,0)$, where $\cdot:[0,1]\times [0,1]\to [0,1]$
is the  product of real numbers and  $([0,1],\oplus,^*,0)$ is the standard MV-algebra.
In the sequel we prove that the variety of Riesz MV-algebras is generated by $[0,1]$, i.e.
 an identity holds in any Riesz MV-algebra if and only if it holds in the standard Riesz MV-algebra
$[0,1]$. Our approach follows closely the proof of  Chang's completeness theorem for \L ukasiewicz logic \cite{chang}.
To any sentence in the first-order theory of Riesz MV-algebras we associate a sentence in the first-order theory of Riesz spaces such that the satisfiability is preserved by the $\Gamr$ functor. 
The  first-order theory of Riesz MV-algebras, as well as the theory of Riesz spaces, are obtained considering for 
each scalar $r$ an unary function $\rho_r$ which denotes in a particular model the scalar multiplication by $r$, i.e.
$x\stackrel{\rho_r}{\mapsto} rx$. In the following, the language of Riesz MV-algebras is ${\cal L}_{RMV}=\{\oplus,^*,0,\{\rho_r\}_{r\in [0,1]}\}$ and
the   language of Riesz spaces is ${\cal L}_{Riesz}=\{\leq, +,-,\vee,\wedge, 0,\{\rho_r\}_{r\in \mathbb{R}}\}$.

Let $t(v_{1},\ldots,v_{k})$ be  a term of ${\cal L}_{RMV}$ and $v$ a propositional variable different
from $v_{1}$, $\ldots$, $v_{k}$. We define $\ti{t}$ as follows:

- if $t=0$ then $\ti{0}$ is $0$,

- if $t=v$  then $\ti{t}$ is $v$

- if $t=t_{1}^*$ then  $\ti{t}$ is $v-\ti{t_{1}}$,

- if $t=t_{1}\oplus t_{2}$ then $\ti{t}$ is $(t_{1}+t_{2})\wedge v$,

- if $t=\rho_r(t_1)$ then $\ti{t}$ is $\rho_r(\ti{t_1})$.\\

Let $\varphi(v_{1},\ldots,v_{k})$ be  a formula of ${\cal L}_{RMV}$ such that all the free and bound 
variables of $\varphi$ are in $\{v_{1},\ldots,v_{k}\}$ and $v$ a propositional variable 
different from $v_{1}$, $\ldots$, $v_{k}$. We define $\ti{\varphi}$ as follows:

- if $\varphi$ is $t_{1}=t_{2}$ then $\ti{\varphi}$ is $\ti{t_1}=\ti{t_2}$,

- if $\varphi$ is ${\neg}\psi$ then $\ti{\varphi}$ is ${\neg}\ti{\psi}$,

- if $\varphi$ is $\psi{\vee}\chi$ then $\ti{\varphi}$ is $\ti{\psi}{\vee}\ti{\chi}$ and 
similarly for
\vspace*{-0.2cm}
\bce
  ${\wedge}$, ${\sra}$, $\leftrightarrow$,
\ece
\vspace*{-0.2cm}

- if $\varphi$ is $(\forall v_{i})\psi$ then $\ti{\varphi}$ is 
$\forall v_{i}((0\leq v_{i}){\wedge} (v_{i}\leq v){\sra}\ti{\psi})$,

- if $\varphi$ is $\exists v_{i}\psi$ then $\ti{\varphi}$ is 
$\exists v_{i}((0\leq v_{i}){\wedge} (v_{i}\leq v){\sra}\ti{\psi})$.

\vspace*{0.2cm}

\noindent Thus to any formula $\varphi(v_{1},\ldots, v_{k})$ of ${\cal L}_{RMV}$ we associate a formula
$\ti{\varphi}(v_{1},\ldots, v_{k},v)$ of ${\cal L}_{Riesz}$. As a consequence, to any sentence 
$\sigma$ of ${\cal L}_{RMV}$ corresponds a formula with only one free variable $\ti{\sigma}(v)$
of ${\cal L}_{Riesz}$.

\bprop\label{sat} Let $(V,u)$ be a Riesz space with strong unit and $R=\Gamr(V,u)$. If $\sigma$ is a sentence 
in the first-order theory of Riesz MV-algebras then 
\bce
$R\models\sigma$ if and only if $V\models\ti{\sigma}[u]$.
\ece
\eprop
\begin{proof} By structural induction on terms it follows that $t[a_{1},\ldots, a_{n}]=\ti{t}[a_{1},\ldots, a_{n},u]$
whenever  $t(v_{1},\ldots,v_{n})$ is a  term of ${\cal L}_{RMV}$ and $a_{1}$, $\ldots$, $a_{n}\in R$. The rest of the proof is straightforward.  
\end{proof}

\bthm\label{comp}
An  equation $\sigma$ in the theory  of Riesz MV-algebras holds in all Riesz MV-algebras if and only if 
it   holds in the standard Riesz MV-algebra $[0,1]$.
\ethm
 \begin{proof}  One implication is obvious. To prove the other one,  let $R$ be  a Riesz MV-algebra such that $R\not\models\sigma$. Since $R\simeq\Gamr(V,u)$ for some Riesz space with strong unit $(V,u)$, we have that $\Gamr(V,u)\not\models\sigma$. Using Proposition \ref{sat}, we infer that  $V\not\models\ti{\sigma}[u]$ in 
 the theory of Riesz spaces.  Since the order relation in any lattice can be expressed equationally, we note that 
 $\ti{\sigma}(v)$ is a quasi-identity.  By \cite[Corollary 2.6]{vanAlten} a quasi-identity is satisfied by all Riesz spaces if and only if it is satisfied by $\mathbb{R}$. Hence there exists a  real number $c\geq 0$ such that 
 $\mathbb{R}\not\models\ti{\sigma}[c]$. Since $\mathbb{R}\models\ti{\sigma}[0]$, we get $c>0$. If follows that 
 $f:\mathbb{R}\to \mathbb{R}$ defined  by $f(x)\mapsto x/c$ is an automorphism of Riesz spaces. We infer that
 $\mathbb{R}\not\models\ti{\sigma}[1]$, so $[0,1]\not\models\sigma$.  
\end{proof}

\bfact
The variety of Riesz MV-algebras is generated by the standard model $[0,1]$ in the language of MV-algebras enriched  with unary operations $x\mapsto rx$ for any $r\in [0,1]$. This features are reflected by the propositional calculus $\rmvlog$ presented in Section \ref{propc}, whose Lindenbaum-Tarski algebra is a  Riesz MV-algebra.

Since the class of Riesz MV-algebras is  a variety, free structures exist. The free Riesz MV-algebra with $n$ free generators is characterized in Corollary \ref{freeR}.
\efact

\section{Norm-complete Riesz MV-algebras}\label{sec-norm}

Let  $(V,u)$ be  a unital  Riesz space and  define
\bcev
 $\|\cdot\|_u:V\to {\mathbb R}$  by

$\|x\|_u=\inf\{\alpha\geq 0\mid  |x|\leq\alpha u\}$ for any $x\in V$.
\ecev  
Then $\|\cdot\|_u$ is a seminorm \cite[Proposition 1.2.13]{Meyer} and
\bcev
$|x|\leq |y|$ implies $\|x\|_u\leq \|y\|_u $ for any $x$, $y\in V$. 
\ecev

\bfact\label{obs}
If  $(V,u)$ is a unital  Riesz space, then 
\bcev  $\|x\|_u=\inf\{\alpha\in [0,1]\mid  x\leq\alpha u\}$ for any $x\in [0,u]$.
\ecev
\efact

  This fact leads us to the following definition, which was suggested by V. Marra (private communication).

\bdfn\cite{VM}
If $R$ is a Riesz MV-algebra then the {\em unit seminorm} $\|\cdot\|:R\to [0,1]$  is defined by 
\bcev
$\|x\|=\inf\{r\in [0,1]\mid  x\leq r1\}$ for any $x\in R$.
\ecev  
\edfn

\bfact\label{fnorm}
If $R_1$ and $R_2$ are Riesz MV-algebras and \mbox{$f: R_1\to R_2$} is a morphism, then $x\leq r1$ in $R_1$ implies $f(x)\leq r1$ in $R_2$, so 
$\|f(x)\|\leq \|x\|$ for any 
$x\in R_1$.    If $f$ is injective then  $\|f(x)\|=\|x\|$ for any 
$x\in R_1$.

This fact allows us to infer properties of the unit seminorm in Riesz MV-algebras directly from   the properties of the unit seminorm in Riesz spaces.
\efact

\blem\label{norm}
 In any Riesz MV-algebra $R$, the following properties hold for any $x$, $y\in R$ and $r\in [0,1]$.\\
(a) $\|0\|=0$, $\|1\|=1$,\\
(b) $\|x\oplus y\|\leq \|x\|+\|y\|$,\\
(c) $x\leq y$ implies  $\|x\|\leq \|y\|$,\\
(d) $\|rx\|=r\|x\|$,\\
(e)   if $(m-1)x\leq x^* $ then $\|mx\|=m\|x\|$ for any natural number $m\geq 1$.\\ 
\elem
\begin{proof} By Theorem \ref{cat} and Remark \ref{fnorm} we can safely assume that $R$ is  $[0,u]_V$ for some unital Riesz space $(V,u)$. Hence 
(a)-(d) follow from the properties of the unit seminorm in Riesz spaces \cite[25H]{Fremlin}.\\
(e) Note that $(m-1)x\leq x^*$ implies 
\bcev
$\underbrace{x\oplus\cdots\oplus x}_{m}=\underbrace{x+\cdots+ x}_{m}$,
\ecev
 where $+$ is the group addition of $V$, so the desired equality is straightforward.  
\end{proof}

\bex\label{csup}
 If $X$ is a compact Hausdorff space, then $C(X)_u=\{f:X\to [0,1]| f \,\, \mbox{continuous}\}$ is a Riesz MV-algebra  and, 
for any $f\in C(X)_u$, we have 
\bcev
$\|f\|=\inf\{r\in [0,1]\mid f(x)\leq r \,\,\forall x\in X\}=\sup\{f(x)|x\in X\}=\|f\|_\infty$.
\ecev 
\eex

Recall that 
 an {\em M-space} is a unital Riesz space $(V,u)$ that is norm-complete with respect to the unit norm.

\bex
If $X$ is a compact Hausdorff space and \bce $C(X)=\{f:X\to {\mathbb R}\mid f \,\, \mbox{continuous}\}$,\ece
then $(C(X),{\mathbf1})$ is an M-space, where ${\mathbf 1}$ is the constant function 
${\mathbf 1}(x)=1$ for any $x\in X$.
\eex

The above example is fundamental, as proved by Kakutani's representation theorem.

\bthm\label{kakutani}\cite{KaM}
For any M-space $(V,u)$ there exists a compact Hausdorff space $X$ such that $(V,u)$  is isomorphic with  $(C(X),{\mathbf 1})$.
\ethm

Let us denote by ${\mathbb M}{\mathbb U}$ the category of M-spaces with unit-preserving morphisms and by 
${\mathbb K}{\mathbb H}aus{\mathbb S}p$ the category of compact Hausdorff spaces with continuous maps. 

\bthm\label{khau}\cite{Bana,SS}
The category ${\mathbb K}{\mathbb H}aus{\mathbb S}p$ is dual to the category ${\mathbb M}{\mathbb U}$. 
\ethm 

\noindent We characterize in the sequel those Riesz MV-algebras that are, up to isomorphism, unit intervals in M-spaces. 
Note that, on any Riesz MV-algebra $R$, we can define
\bcev
$\delta_{\|\cdot\|}(x,y)=\|d(x,y)\|$ for any $x$, $y\in R$.
\ecev
By Lemmas \ref{dist} and \ref{norm} it follows that $\delta_{\|\cdot\|}$ is a pseudometric on $R$. 

\bdfn
We say that a Riesz MV-algebra $R$  is {\em norm-complete} if $(R,\delta_{\|\cdot\|})$ is a complete metric space.
\edfn

\bthm\label{rnorm} If $(V,u)$ is a unital Riesz space  then the following are equivalent: \\
(i)  $(V,u)$ is an M-space,\\ 
(ii)  $\Gamr(V,u)$ is a norm-complete Riesz MV-algebra.

\ethm
\begin{proof} We denote  $R=\Gamr(V,u)$.\\
(i)$\Rightarrow$(ii) By  Remark \ref{obs}, $\|x\| = \|x\|_{u}$ for any $x\in [0,u]$.
In consequence, any Cauchy sequence w.r.t $\|\cdot\|$ from $R\Gamr(V,u)$ is a Cauchy sequence w.r.t  $\|\cdot\|_{u}$ in $(V,u)$ and we use the fact that $(V,u)$ is norm-complete. \\
(ii)$\Rightarrow$(i)  Let $(v_n)_n$ be a Cauchy sequence in $V$ w.r.t.
 $\|\cdot\|_{u}$ such that $v_n\geq 0$ for any $n$. It follows that it is bounded, i.e. there is $y\in V$  and 
$\|v_n\|_{u}\leq \|y\|_{u}$ for any $n$. We get $v_n\leq \|v_n\|_{u} u\leq \|y\|_{u}u$ for any $n$, so there exists a natural number  $k$ such that $v_n\leq ku$ for any $n$. By Lemma \ref{l:goodsequence},
\bcev
$v_n=v_{n_1}+\cdots+v_{n_k}$,  where $v_{n_i}=((v_n-(i-1)u)\vee 0)\wedge u\in [0,u]$
\ecev
for any $i\in\{1,\ldots,k\}$. Since $v_{n_i}\in [0,u]$ we get 
\bce $ \|v_{n_i}\|_{u}=\|v_{n_i}\|$ for any $n$ and $i\in \{1,\ldots,k\}$.\ece
 One can easily see that 
$(v_{n_i})_n$ is  a Cauchy sequence in $R=[0,u]$ for any $i\in \{1,\ldots,k\}$. Since $R$ is norm-complete it follows that, for 
any $i\in \{1,\ldots,k\}$ there is $w_i\in R$ such that $\lim_n d(v_{n_i},w_i)=\lim_n\|v_{n_i}-w_i\|_u=0$.  If $w=w_1+\cdots +w_k$ then 
\bcev
$\|v_n-w\|_u\leq\|v_{n_1}-w_1\|_u+\cdots +\|v_{n_k}-w_k\|_u$   
\ecev
for any  $n$, so $\lim_n\|v_n-w\|=0$ and $(v_n)_n$ is convergent w.r.t. $\|\cdot\|_u$ in $V$. 

Recall that $v=(v\vee 0) -((-v)\vee 0)$  for any $v\in V$, so the convergence of arbitrary Cauchy sequences reduces to 
the convergence of positive Cauchy sequences.  
\end{proof}

Denote by ${\mathbb U}{\mathbb R}{\mathbb M}{\mathbb V}$ the category of norm-complete Riesz MV-algebras, which is a full subcategory of  ${\mathbb R}{\mathbb M}{\mathbb V}$.  
By Remark \ref{fnorm}, the norm-preserving morphisms coincide with the monomorphisms of  ${\mathbb U}{\mathbb R}{\mathbb M}{\mathbb V}$.

Using Theorem \ref{rnorm}, the functor $\Gamr$  yields the following categorical equivalence.

\bcor
The categories ${\mathbb U}{\mathbb R}{\mathbb M}{\mathbb V}$ and ${\mathbb M}{\mathbb U}$ are equivalent.
\ecor

\bcor
The categories ${\mathbb U}{\mathbb R}{\mathbb M}{\mathbb V}$ and ${\mathbb K}{\mathbb H}aus{\mathbb S}p$ are dually  equivalent.
\ecor

\bfact
Following \cite[Chapter IV]{SS} and \cite{Bana}, the functors establishing the above equivalences are defined on objects as follows:
\bcev
$R\mapsto Max(R)$ and $X\mapsto C(X)_u$
\ecev
 for any  norm-complete Riesz space $R$ and compact Hausdorff space $X$, where $Max(R)$ is the set of all maximal ideals of $R$. 
\efact

In the sequel, we connect our result with the Gelfand-Naimark duality for C$^*$-algebras \cite{GN}. Recall that  MV-algebras are related with AF C$^*$-algebras in \cite{Meq}, but in this case the K-theory is used. 

Denote ${\mathbb C}^*$ the category whose objects are commutative unital C$^*$-algebras and whose morphisms are unital 
C$^*$-algebra morphisms.  

\bthm\cite[Chapter 1.1]{MK}
The categories  ${\mathbb C}^*$  and ${\mathbb K}{\mathbb H}aus{\mathbb S}p$ are dually  equivalent.
\ethm

As a corollary we infer immediately that the categories of commutative unital  C$^*$-algebras and norm-complete Riesz MV-algebras are equivalent. 

\bcor
The categories  ${\mathbb C}^*$  and  ${\mathbb U}{\mathbb R}{\mathbb M}{\mathbb V}$ are  equivalent.
\ecor

\section{The  propositional calculus  ${\mathbb R}{\cal L}$}\label{propc}
\noindent We denote by  ${\cal L}_\infty$ the $\infty$-valued propositional \L ukasiewicz logic. Recall that   ${\cal L}_\infty$  has $\neg$ (unary) and $\sra$ (binary) as  primitive  connectives and, for any $\varphi$ and $\psi$ we have:

$\vp\vee \psi\defi (\vp\sra\psi)\sra\psi$, $\vp\wedge\psi\defi\neg(\neg\vp\vee\neg\psi)$,

$\vp\ba\psi\defi (\vp\sra\psi)\wedge (\psi\sra\vp)$.\\
The language of ${\mathbb R}{\cal L}$ contains the language of  ${\cal L}_\infty$ 
and a family of unary connectives   $\{\nabla_r| r\in [0,1]\}$. 
We denote by $Form({\mathbb R}{\cal L})$ the set of formulas, which are defined inductively as usual.

\bdfn\label{7.2} An {\em axiom} of ${\mathbb R}{\cal L}$ is
  any formula that is an axiom of ${\cal L}_\infty$ and any formula that has one of the following forms,  where $\vp$, $\psi$, $\chi\in Form({\mathbb R}{\cal L})$ and
\mbox{$r$, $q\in [0,1]$:}

\begin{enumerate}[(RL1)]
\item $\nabla_{r}(\vp\sra\psi)\ba (\nabla_{r}\vp\sra\nabla_{r}\psi)$; 
\item $\nabla_{(r\odot q^{*})}\vp\ba (\nabla_{q}\vp\sra\nabla_{r}\vp)$;
 \item $\nabla_{r}\nabla_{q}\vp\ba \nabla_{(rq)}\vp$; 
\item $\nabla_{1}\vp\ba\vp$,
\end{enumerate}

\noindent The  deduction rule of ${\mathbb R}{\cal L}$ is
 {\em modus ponens} and provability is   defined as usual. 
\edfn

\bfact  $\{\vp\}\vdash \nabla_r\vp$ is a derived deduction rule for any
 $r\in [0,1]$.
\efact

We recall the usual construction of the Lindenbaum-Tarski algebra.
The equivalence relation $\equiv$ is defined on $Form({\mathbb R}{\cal L})$ as follows:
\begin{center}
$\vp\equiv \psi$ iff $\vdash\vp\sra\psi$ and
$\vdash \psi\sra\vp$.
\end{center}
\noindent
We denote by  $[\vp]$ the equivalence class of a formula $\vp$ and we define
on $Form({\mathbb R}{\cal L})/{\equiv}$  the following operations:\\
$[\vp]^{*}= [\neg\vp]$,\\
$[\vp]\oplus[\psi]=[\neg\vp\sra\psi]$, $[\vp]\odot[\psi]=[\neg(\vp\sra\neg \psi)]$,\\
$0=[\neg(v_1\sra v_1)]$, $1=0^*=[v_1\sra v_1]$. \\
In order to define the scalar multiplication we introduce new  connectives:
\bce
$\Delta_r\vp \defi \neg(\nabla_r\neg\vp)$  

\ece
and we set $r [\varphi]= [\Delta_r\vp]$   for any $r\in [0,1]$ and $\vp$ formula of ${\mathbb R}{\cal L}$.

\medskip

\bprop
The Lindenbaum-Tarski algebra \bce  $RL=(Form({\mathbb R}{\cal L})/{\equiv},\cdot, \oplus,^{*}, 0)$\ece
 is a Riesz MV-algebra.
\eprop
\begin{proof}
The axioms  (RL1)-(RL4) are logical expressions of the  duals of (RMV1)-(RMV4).   We prove in detail that $RL$ satisfies (RMV1). If  $\varphi$ and $\psi$ are two formulas and $r\in [0,1]$ then, by (RL1), we get
\bce $[\nabla_r(\neg\vp\sra\neg\psi)]=[\nabla_r\neg\vp\sra\nabla_r\neg\psi]$.\ece It follows that:
\bce
$[\nabla_r\neg(\neg\vp\odot\psi)]=[\neg (\nabla_r\neg\vp\odot\neg \nabla_r\neg\psi)]$

$[\neg \nabla_r\neg(\neg\vp\odot\psi)]=[\nabla_r\neg\vp\odot\neg \nabla_r\neg\psi]$

$[\Delta_r(\neg\vp\odot\psi)]=[\neg\Delta_r\vp\odot\nabla_r\psi]$

$r([\vp]^*\odot[\psi])=(r[\vp])^*\odot (r[\psi])$,
\ece
so (RMV1) holds in $RL$.  
\end{proof}
\medskip

Let  $R$ be an Riesz MV-algebra. An {\em evaluation} is a function
\mbox{$e:Form({\mathbb R}{\cal L})\sra R$} which satisfies the following conditions for any $\vp$, $\psi\in Form({\mathbb R}{\cal L})$ and $r\in [0,1]$:\\
(e1) $e(\vp\sra\psi)=e(\vp)^*\oplus  e(\psi)$,\\
(e2) $e(\neg\vp)=e(\vp)^{*}$,\\
(e3) $e(\nabla_{r}\vp)=(re(\vp)^*)^{*}$.

\medskip

As a consequence of Theorem \ref{comp}, the propositional calculus $\rmvlog$ is complete with respect to $[0,1]$.

\bthm\label{thcomp}
For a formula $\vp$ of  $\rmvlog$ the following are equivalent:\\
(i) $\vp$ is provable in $\rmvlog$,\\
(ii) $e(\vp)=1$  for any Riesz MV-algebra $R$ and for any evaluation $e:Form({\mathbb R}{\cal L})\sra R$,\\
\mbox{(iii) $e(\vp)=1$  for any evaluation $e:Form({\mathbb R}{\cal L})\sra [0,1]$.}

 \ethm

\bfact
The system  ${\mathbb R}{\cal L}$ is a conservative extension of ${\cal L}_\infty$, i.e. a formula $\vp$ of ${\cal L}_\infty$ is a theorem of ${\cal L}_\infty$ if and only if it is a theorem of ${\mathbb R}{\cal L}$ . Since any proof in ${\cal L}_\infty$  is also a proof 
in ${\mathbb R}{\cal L}$, one implication is obvious. To prove the other one, assume that  $\vp$ is a formula of ${\cal L}_\infty$ which is not a theorem of ${\cal L}_\infty$. Hence there exists an evaluation $e^\prime: Form({\cal L}_\infty)\to [0,1]$ such that 
$e^\prime(\vp)\neq 1$.  Let  $e:Form({\mathbb R}{\cal L})\sra [0,1]$ the unique evaluation in ${\mathbb R}{\cal L}$ such that 
$e(v)=e^\prime(v)$ for any propositional variable $v$. By structural induction on formulas one can prove that 
$e(\psi)=e^\prime(\psi)$ for any $\psi\in Form({\cal L}_\infty)$. It follows that  $e(\vp)=e^\prime(\vp)\neq 1$, so 
$\vp$ is not a theorem of   ${\mathbb R}{\cal L}$. 
\efact

\bfact
A formula  $\varphi$  with variables  from $\{v_1,\ldots,v_n\}$ uniquely defines a  {\em term function}:
\bce
$\fct{\vp}:[0,1]^n\to [0,1]$, $\fct{\vp}(x_1,\ldots,x_n)=e(\vp)$,
\ece
 where $e$ is an evaluation such that $e(v_i)=x_i$ for any $i\in\{1,\ldots,n\}$.
 By Theorem \ref{thcomp} it follows that $[\vp]=[\psi]$ if and only if $\fct{\vp} = \fct{\psi}$.
\efact

\section{Term functions and piecewise linear functions}\label{normal}

 In the following, we characterize the class of functions that can be defined by formulas in ${\mathbb R}{\cal L}$.

\bdfn
Let $n>1$ be a natural number.
A {\em piecewise linear function}  is a function  \mbox{$f:{\mathbb R}^n\to {\mathbb R}$}  for which 
  there exists a finite number of affine functions   \bce$q_1$, $\ldots$, $q_k:\real^n\to \real$\ece and
for any $(x_1,\ldots, x_n)\in {\mathbb R}^n$  there is $i\in \{1,\ldots, k\}$ such that
$f(x_1,\ldots, x_n)=q_i(x_1,\ldots, x_n)$. We say that $q_1$, $\ldots$, $q_k$ are the {\em components} of $f$.

We denote by$PL_n$ the set of  all  continuous functions $f:[0,1]^n\to [0,1]$ that are  piecewise linear.

For the rest of the paper, all  piecewise linear functions are continuous.  
  
\edfn

\bthm\label{mcn}
If   $\vp$ is a formula  of  ${\mathbb R}{\cal L}$ with propositional  variables from    $\{v_1,\cdots,v_n\}$ then
$\fct{\varphi}\in PL_n$.
\ethm
\begin{proof}  We prove the result by structural induction on formulas.

If $\varphi$ is $v_i$ for some $i\in\{1,\ldots ,n\}$ then 
$\fct{\varphi}=\pi_i$ (the $i$-th projection).

If $\varphi$ is $\neg\psi$ and  $q_{1}$, $\ldots$, $q_{s}$  are  the components 
 of $\fct{\psi}$, then  $1-q_{1}$, $\ldots$, $1-q_{s}$ are the components of $\fct{\varphi}$.
 
Assume $\varphi$ is $\psi\sra\chi$. If  $q_{1}$, $\ldots$, $q_{m}$  are the  
components of $\fct{\psi}$ and 
$p_{1}$, $\ldots$, $p_{k}$ are the components of $\fct{\chi}$, then $\fct{\varphi}$ is defined by 
$\{1\}\cup \{s_{ij}\}_{i,j}$, where $s_{ij}=1-q_i+p_j$ for any $i\in\{1,\ldots ,s\}$ and $j\in\{1,\ldots ,r\}$.

 If $\varphi$ is   $\Delta_r\psi$ for some $r\in [0,1]$ and  $q_{1}$, $\ldots$, $q_{s}$  are  the components
 of $\fct{\psi}$, then  $1-r+rq_{1}$, $\ldots$, $1-r+rq_{s}$  are the components of $\fct{\varphi}$.  
\end{proof}

\bfact
The  continuous piecewise linear functions \\$f:[0,1]^n\ra [0,1]$ with {\em integer coefficients} are called {\em McNaughton functions} and they
are in one-one correspondence with the  formulas of \L ukasiewicz logic  by McNaughton theorem \cite{McN}.  The continuous piecewise linear functions with {\em rational coefficients} correspond to formulas of Rational \L ukasiewicz logic, a propositional calculus developed in \cite{Brunella} that has divisible MV-algebras as models. 
In Theorem \ref{main} we prove that any continuous piecewise linear function  with  {\em real coefficients} $f:[0,1]^n\to [0,1]$ is the term function of a  formula from ${\mathbb R}{\cal L}$.
\efact

For now on we define $\sq :{\mathbb R}\to [0,1]$ by
 \begin{center}
$\sq(x)=(x\vee 0)\wedge 1$ for any $x\in {\mathbb R}$.
\end{center}

\medskip

\blem\label{calc}
For any $x$, $y\in {\mathbb R}$ the following hold:\\
(a) $(x\vee 0) +(y\vee 0)\geq (x+y)\vee 0$,\\
(b)  $x\geq 0$ iff $\sq(-x)=0$,\\
(c) $\sq(x)=\sq(x\vee 0)$.
\elem
\begin{proof}
(a) $(x\vee 0) +(y\vee 0)=(x+y)\vee x\vee y\vee 0\geq (x+y)\vee 0$.\\
(b) $\sq(-x)=0$ iff $((-x)\vee 0)\wedge 1=0$ iff $(-x)\vee 0=0$ iff $-x\leq 0$ iff $x\geq 0$.\\
(c) $\sq(x\vee 0)=(x\vee 0\vee 0)\wedge 1=(x\vee 0)\wedge 1 =\sq(x)$.  
\end{proof}

In the following we generalize some results from \cite[Lemma 3.1.9]{CDOM}.

\blem\label{ppmain} If $g:[0,1]^n\to {\mathbb R}$  and $h:[0,1]^n\to [0,1]$ then the following properties hold.\\
(a) $\sq\circ (g+h) = ((\sq\circ g)\oplus h)\odot (\sq\circ (g+1))$.\\
(b)   $\sq\circ (1-g)=1-(\sq\circ g)$
\elem
\begin{proof}
Let   $\bm{x}=(x_1,\cdots, x_n)$ be an element from  $ [0,1]^n$.\\
(a) If $g(\bm{x})>1$  then  $g(\bm{x})+1>1$ and $g(\bm{x})+h(\bm{x})>1$. It follows that
\bce
$\sq({g}(\bm{x}))=\sq({(g+1)}(\bm{x}))=\sq((g+h)(\bm{x}))=1$,
\ece
\noindent so the intended identity is obvious.

If $g(\bm{x})\in[0,1]$ then
$\sq({g}(\bm{x}))=g(\bm{x})$ and $\sq({(g+1)}(\bm{x}))=1$   for any $\bm{x}\in [0,1]^n$, so

\begin{center}$\begin{array}{ll}
\sq({(g+h)}(\bm{x}))&= h(\bm{x})\oplus g(\bm{x})\\
    &=h(\bm{x})\oplus \sq({g}(\bm{x}))=((\sq\circ{g})\oplus h)(\bm{x})\\
    &=(((\sq\circ{g})\oplus h)\odot 1)(\bm{x})\\
    &=(((\sq\circ{g})\oplus h)\odot (\sq\circ{(g+1)}))(\bm{x}).
\end{array}$\end{center}

Assume that $g(\bm{x})<0$, so $\sq({g}(\bm{x}))=0$.   We have to prove that
$\sq({(g+h)}(\bm{x}))=h(\bm{x})\odot \sq({(g+1)}(\bm{x}))$.

If  $g(\bm{x})\leq -1$ then $\sq(g(\bm{x}))=\sq({(g+1)}(\bm{x}))=0$ and $g(\bm{x})+h(\bm{x})\leq -1+h(\bm{x})\leq 0$,
so
\bce
 $\sq({(g+h)}(\bm{x}))=0=h(\bm{x})\odot 0=(h\odot (\sq\circ{(g+1)}))(\bm{x})$.
\ece

If $g(\bm{x})\in (-1,0)$ then  $\sq({(g+1)}(\bm{x}))=(g+1)(\bm{x})=g(\bm{x})+1$, so we get

\begin{center}
$\begin{array}{ll}
(\sq\circ{(g+h)})(\bm{x})&= 0\vee (1\wedge (h(\bm{x})+g(\bm{x})))\\ & =0\vee  (h(\bm{x})+g(\bm{x}))\\
     &=0\vee (h(\bm{x})+g(\bm{x})+1-1)\\ & =h(\bm{x})\odot (g(\bm{x})+1)\\
    &=h(\bm{x})\odot \sq({(g+1)}(\bm{x}))\\ & = (h\odot (\sq\circ{(g+1)}))(\bm{x}).
 \end{array}$
\end{center}

\noindent  (b)  If $g(\bm{x})<0$ then $\sq({g}(\bm{x}))=0$ and $$\sq({(1-g)}(\bm{x}))=1=1-0=1-\sq({g}(\bm{x})).$$

If $g(\bm{x})\in [0,1]$ then $(1-g)(\bm{x})\in [0,1]$, so  $$\sq({(1-g)}(\bm{x}))=(1-g)(\bm{x}) =1-g(\bm{x})=1-\sq({g}(\bm{x})).$$

If  $g(\bm{x})>1$ then $\sq({g}(\bm{x}))=1$ and \bce$\sq({(1-g)}(\bm{x}))=0=1-1=1-\sq({g}(\bm{x}))$.  \ece
\end{proof}

\bprop\label{pmain}
For   any affine  function  $f:[0,1]^n\to {\mathbb R}$  there exists a formula $\varphi$ of  ${\mathbb R}{\cal L}$ such that $\sq\circ f=\fct{\vp}$.
\eprop
\begin{proof}
Let $f:[0,1]^n\to {\mathbb R}$  be an affine function, i.e. there are $c_0$, $\ldots$, $c_n\in {\mathbb R}$ such that
 \[f(x_1,\ldots, x_n)=c_nx_n+\cdots+c_1x_1+c_0 \]
 for any  $(x_1,\ldots, x_n)\in [0,1]^n $.
 Note that for  any  $c\in {\mathbb R}$  there is a natural number $m$ such that
 $c=r_1+\cdots+r_m$  where  $r_1$, $\ldots$, $r_m\in [-1,1]$. Hence we  assume that
\[ f(x_1,\ldots, x_n)=r_my_m+\cdots r_{p+1}y_{p+1}+r_{p}+\cdots+r_1\]
where $m\geq 1$ and $ 0\leq p\leq m$ are natural numbers,   $r_j\in [-1,1]\setminus\{0\}$ for any $j\in\{1,\ldots, m\}$ and\\
 $y_j\in\{x_1,\ldots,x_n\}$ for any $j\in\{p+1,\cdots, m\}$.
\medskip

We prove the theorem by induction on $m\geq 1$.
Let us denote $\bm{x}=(x_1,\cdots, x_n)$  an element from  $ [0,1]^n$
\medskip

\noindent {\em Initial step} $m=1$. We have
$f(\bm{x})=r$ for any $\bm{x}\in [0,1]^n$ or  $f(\bm{x})=rx_i$ for any $\bm{x}\in [0,1]^n$,
 where $r\in  [-1,1]\setminus\{0\}$ and  $i\in\{1,\ldots,n\}$.
If  $r\in[-1,0)$ then $\sq\circ f=0$  so $\sq\circ f=\fct{\varphi}$ for $\vp=v_1\odot \neg v_1$ .
If $r\in (0,1]$ then  $f=\sq\circ f$.  It follows that $f=\fct{\vp}$ where $\vp=\nabla_r(v_1\sra v_1)$
if $f(\bm{x})=r$ for any $\bm{x}\in [0,1]^n$ and  $\vp={\nabla_r v_i}$
if   $f(\bm{x})=rx_i$  for any $\bm{x}\in [0,1]^n$.

\medskip

\noindent{\em Induction step}. We take  $f=g+h$ where $\sq\circ g=\fct{\varphi}$ for some formula $\varphi$ and there are
 $r\in  [-1,1]\setminus\{0\}$ and  $i\in\{1,\ldots,n\}$ such that
$h(\bm{x})=r$  for any $\bm{x}\in [0,1]^n$, or  $h(\bm{x})=rx_i$ for any $\bm{x}\in [0,1]^n$.
We consider two cases.
\medskip

\noindent{\em Case 1.} If   $r\in (0,1]$ then $h:[0,1]^n\to [0,1]$ so
\bce $\sq\circ f =((\sq\circ g)\oplus h)\odot (\sq\circ (1+g))$\ece  by Lemma \ref{ppmain} (a).
  Following the initial step, there is a formula $\psi$ such that $h=\fct{\psi}$. Note that  $1+g=1-(-g)$ and,  since   the induction hypothesis holds for $(-g)$,
there is a  formula $\chi$ such that $\sq\circ (-g)=\fct{\chi}$. In consequence, by Lemma \ref{ppmain} (b),
$\sq\circ (1+g)=1-\fct{\chi}=\fct{\neg\chi}$.
We get $\sq\circ f=\fct{\theta}$ where $\theta={(\vp\oplus\psi)\odot\neg\chi}$.
\medskip

\noindent{\em Case 2.} If $r\in [-1,0)$, then   $g+h=(g-1)+(1+h)$ and $1+h:[0,1]^n\to [0,1]$. By  Lemma \ref{ppmain} (a) we get
\bce
$\sq\circ f =((\sq\circ (g-1))\oplus (1+h))\odot (\sq\circ g)$.
\ece
\noindent Following the initial step, there is a formula $\psi$ such that
\bce $-h=\fct{\psi}$, so $1+h=1-(-h)=\fct{\neg\psi}$.\ece
In the sequel we have to find a formula $\chi$ that corresponds to $\sq\circ (g-1)$, where
\[ g(\bm{x})=r_my_m+\cdots r_{p+1}y_{p+1}+r_{p}+\cdots+r_1\]
with  $r_j\in [-1,1]\setminus\{0\}$ for any $j\in\{1,\cdots, m\}$ and  $y_j\in\{x_1,\ldots,x_n\}$ for any $j\in\{p+1,\cdots, m\}$.

\medskip

\noindent{\em Case 2.1.} If  $r_j\leq 0$ for any $j\in\{1,\cdots, m\}$ then \\$g-1\leq 0$, so $\sq\circ (g-1)=0= \fct{\chi}$ with
 $\chi=v_1\odot \neg v_1$.
\medskip

\noindent{\em Case 2.2.} If there is $j_0\in\{1,\cdots, p\}$ such that $r_{j_0}>0$, then it follows that
\bce
$(g-1)(\bm{x})=r_my_m+\cdots r_{p+1}y_{p+1}+r_{p}+\cdots+(r_{j_0}-1)+\cdots +r_1$
\ece
and $r_{j_0}-1\in [-1,0)$, so the induction hypothesis applies to $g-1$.  In consequence,  there exists a formula $\chi$ such that
$\sq\circ (g-1)=\fct{\chi}$.
\vspace*{0.2cm}

\noindent{\em Case 2.3.} If there is $j_0\in\{p+1,\cdots, m\}$ such that $r_{j_0}>0$, then we set  $h_0(\bm{x})=r_{j_0}y_{j_0}$ and
$$g_0(\bm{x})=g(\bm{x})-r_{j_0}y_{j_0}-1.$$
It  follows that $g-1=g_0+h_0$ such that $g_0$ satisfies the induction hypothesis and $h_0:[0,1]^n\to [0,1]$. We are in the hypothesis of {\em Case 1}, so
there exists a formula $\chi$ such that $\sq\circ(g-1)=\fct{\chi}$.
\vspace*{0.2cm}

Summing up, we get  $\sq\circ (g+h)  =\fct{\theta}$ with $\theta={((\chi\oplus\neg\psi)\odot\vp)}$.  
\end{proof}

\bthm\label{main}
 For any  $f:[0,1]^n\to [0,1]$ from $PL_n$  there is a formula $\varphi$ \mbox{of  ${\mathbb R}{\cal L}$ such that $f=\fct{\vp}$.}
\ethm
\begin{proof}
Let $f:[0,1]^n\to [0,1]$ be in $PL_n$.
 Using  the Max-Min representation  from \cite{minmax}, there are finite sets $I$ and $J$ such that
$$f=\bigvee_{i\in I}\bigwedge_{j\in J} f_{ij},$$
\noindent where $f_{ij}:[0,1]^n\to{\mathbb R}$ are affine functions. We note that

$$f=\sq\circ f=\bigvee_{i\in I}\bigwedge_{j\in J}( \sq\circ f_{ij}).$$

 By Proposition \ref{pmain},  for any $i\in I$ and $j\in J$ there is a formula $\varphi_{ij}$ such that
$\sq\circ f_{ij}=\fct{\varphi_{ij}}$. In consequence, if we set $\varphi=\bigvee_{i\in I}\bigwedge_{j\in J} \varphi_{ij}$ then
$f=\fct{\vp}$.  
\end{proof}

For any  $n\geq 1$, the set $PL_n$  is a  Riesz MV-algebra with the operations defined  componentwise. 
If  $RL_n$ is the Lindenbaum-Tarski algebra of ${\mathbb R}{\cal L}$ defined on formulas with  variables from $\{v_1,\ldots, v_n\}$, then $RL_n$ is the free Riesz MV-algebra with $n$ free generators  by  standard results in universal algebra \cite{BS}. 
Since  the function $[\varphi]\mapsto\fct{\vp}$ is obviously an  isomorphism between $RL_n$ and $PL_n$ the following corollary is straightforward.

\bcor\label{freeR}
$PL_n$  is the free Riesz MV-algebra with $n$ free generators.
\ecor

\section{Linear combinations of formulas and de Finetti's coherence criterion}\label{linear}

We recall  in the beginning de Finetti's  coherence criterion for boolean events. If $S=\{\varphi_1,\ldots, \varphi_k\}$ is a  set of classical events then a {\em book} is a set
\bce  $\{(\varphi_i,r_i)\mid i\in \{1,\ldots, k\}\}$,
\ece
 where $r_i\in [0,1]$  is a "betting odd"  assigned by a bookmaker for 
$\varphi_i$ for any $i\in \{1,\ldots, k\}$.  The book is coherent if there is no system of bets $\{c_1,\cdots, c_k\}$ which 
causes the bookmaker a sure loss. This means that for any   real numbers  $\{c_1,\cdots, c_k\}$  there exists an evaluation $e:S\to \{0,1\}$ such that $\sum_{i=1}^kc_i(r_i-e(\varphi_i))\geq 0$. De Finetti's coherence criterion \cite{deF} states that the book 
$\{(\varphi_i,r_i)|i\in \{1,\ldots, k\}\}$ is coherent if there is a boolean probability $\mu$ defined on the algebra of events generated by $S$ such that $\mu(\varphi_i)=r_i$ for any $i\in \{1,\ldots,k\}$.  When the underlying logic is \L ukasiewicz logic \cite{MuBook}, the 
events belongs to an MV-algebra and they are evaluated in $[0,1]$. Consequently, the coherence criterion uses states instead of boolean probabilities.

We recall the MV-algebraic approach to de Finetti's notion of coherence \cite{KM} and we prove a similar coherence criterion for Riesz MV-algebras. We also provide a logical expression of the coherence criterion   when the events are represented by formulas  ${\mathbb R}{\cal L}$. In order to  accomplish this task, we initiate the study of linear combinations of formulas in ${\mathbb R}{\cal L}$.

The next result characterizes the {\em states} defined on Riesz MV-algebras. 

\blem
If $R$ is a Riesz MV-algebra then any state  $s: U(R)\to [0,1]$ is also homogeneous:\\
(s3)  $s(r\cdot x)=rs(x)$ for any $r\in [0,1]$, $x\in R$.
\elem
\begin{proof}
We can safely assume that $R=[0,u]_V$ for some unital Riesz space $(V,u)$.
By \cite[Theorem 2.4]{M1}, there is a state $s^\prime :V\to {\mathbb R}$ such that $s^\prime (x)=s(x)$ for any $x\in [0,u]$. 

If $r\in [0,1]\cap{\mathbb Q}$ then $r=\frac{m}{n}$ and $rx=y$ in $[0,u]$  implies $mx=ny$ in $V$. It follows that 
$s^\prime (mx)=s^\prime (ny)$, so $ms(x)=ns(y)$ and we get $s(rx)=s(y)=rs(x)$. 

If $r\in (0,1)$ there are rational sequences $(r_n)_n$ and $(q_n)_n$   such that $r_n\uparrow r$ and $q_n\downarrow r$. Hence 
\bce
$r_ns(x)=s(r_nx)\leq s(rx)\leq s(q_nx)=q_ns(x)$  
\ece
for any $n$ and $x\in [0,1]$. The intended result follows by an application of Stolz-Ces\`{a}ro theorem.  
 
\end{proof}

\bdfn
If $R$ is a Riesz MV-algebra,  a {\em state} on $R$ is a function $s:R\to [0,1]$ which satisfies the conditions (s1), (s2) and (s3)
(additivity, normalization and homogeneity).  
 
The previous lemma asserts that the states of a Riesz  MV-algebra $R$ coincide with the states of its MV-algebra reduct  $\U (R)$. 
\edfn

\medskip

The following definition generalizes de Finetti's notion of coherence  and provides an algebraic approach within \L ukasiewicz logic \cite{KM}.

\bdfn\cite{KM}
If $A$ is an MV-algebra and \mbox{$x_1, \ldots, x_k$} are in $A$  then a map  $\beta:\{x_1,\ldots, x_n\} \to [0,1]$ is {\em coherent}
 if for any $c_1$, $\cdots$, $c_k\in {\mathbb R}$ there exists a morphism of MV-algebras
$e:A\to [0,1]$ such that \bce $\sum_{i=1}^kc_i(\beta(x_i)-e(x_i))\geq 0$.\ece
\edfn

\bthm\label{th:dfmv}\cite[Theorem 3.2]{KM}
If $A$ is an MV-algebra,  $x_1$, $\ldots$, $x_k\in A$ and   $\beta:\{x_1,\ldots, x_n\} \to [0,1]$ then the following are equivalent:\\
(i)  the map $\beta$ is {coherent},\\
(ii) there exists a state $s: A\to [0,1]$ such that
\bce  $s(x_i)=\beta(x_i)$ for any $i\in\{1,\ldots,k\}$,\ece
(iii) there exists $e_1$, $\ldots$, $e_m:A\to [0,1]$ morphisms of MV-algebras such that   $m\leq k + 1$  and
$\beta$ is the restriction of a convex combination of $\{e_1, \ldots, e_m\}$.
\ethm

\medskip

By Remark \ref{fact-mor}, any morphsim of Riesz MV-algebras is just a morphism between the MV-algebra reducts of its domain and codomain. Therefore,
the notion of {\em coherent map} remains unchanged on Riesz MV-algebras and an analogue of Theorem \ref{th:dfmv} can be proved for Riesz MV-algebras as well.

\bcor\label{th:dfrmv}
If $R$ is a Riesz  MV-algebra,  $x_1$, $\ldots$, $x_k\in R$ and   $\beta:\{x_1,\ldots, x_k\} \to [0,1]$ then the following are equivalent:\\
(i)  the map $\beta$ is {coherent},\\
(ii) there exists a state $s: R\to [0,1]$ such that
\bce  $s(x_i)=\beta(x_i)$ for any $i\in\{1,\ldots,k\}$,\ece
(iii) there exists $e_1$, $\ldots$, $e_m:R\to [0,1]$ morphisms such that   $m\leq k + 1$  and
$\beta$ is the restriction of a convex combination of $\{e_1, \ldots, e_m\}$.
\ecor
\begin{proof}
 (i)$\Leftrightarrow$(iii)  and  (ii)$\Rightarrow$(i) follow by Theorem \ref{th:dfmv} applied to the MV-algebra reduct of $R$ and by Remark \ref{fact-mor}.\\
(iii)$\Rightarrow$(ii) There are $\alpha_1$, $\ldots$, $\alpha_m\in [0,1]$ such that
\bce
 $\alpha_1+\cdots + \alpha_m=1$ and

$\beta(x_i)=\alpha_1e_1(x_i)+\cdots + \alpha_me_m(x_i)$ 
\ece
for any  $i\in\{1,\ldots,k\}$. We set $s= \alpha_1e_1+\cdots + \alpha_me_m$  satisfies (s1), (s2) and (s3), so $s:R\to [0,1]$ is the required state.  
\end{proof}
\medskip

The above result is an algebraic version of de Finetti's coherence criterion. In the sequel we provide a logical approach within $\rmvlog$.

We firstly recall de Finetti's coherence criterion for \L ukasiewicz logic ${\cal L}_\infty$.

\bthm\cite{MuBook}
If  $\vp_1, \ldots, \vp_k$ are formulas of  ${\cal L}_\infty$ with variables
  ${v_1,\ldots,v_n}$ and   $r_1, \ldots, r_k\in [0,1]$ then
the following are equivalent:\\
(i)  the book $\{(\vp_i,r_i)\mid i\in\{1,\ldots,k\}\}$ is {coherent},\\
(ii) there exists a state $s:L_n\to [0,1]$ such that
\bce  $s([\vp_i])=r_i$ for any $i\in\{1,\ldots,k\}$,
\ece
where $L_n$ is the Lindenbaum-Tarski algebra of the formulas in $n$ variables and $[\vp]$ is the equivalence class
of the formula $\vp$ in $L_n$.
 \ethm

When we consider ${\mathbb R}{\cal L}$ instead of ${\cal L}_\infty$, we get the following.

\bdfn\label{book}
If  $\vp_1,\ldots,\vp_k$ are formulas of  $\rmvlog$  and  \\ 
$r_1$, $\ldots$, $r_k\in [0,1]$ then
the book \mbox{$\{(\vp_i,r_i)\mid i\in\{1,\ldots,k\}\}$} is {\em coherent} if for any $c_1$, $\cdots$, $c_k\in {\mathbb R}$ there exists an evaluation
$e:Form(\rmvlog )\to [0,1]$ such that
\bce $\sum_{i=1}^kc_i(r_i-e(\vp_i))\geq 0$.
\ece
\edfn

In order to characterize  coherence in logical terms  within ${\mathbb R}{\cal L}$ we introduce the notion of {\em quasi-linear combination} of piecewise linear functions.

\medskip

In the sequel we assume that  $n$ is a natural number and all the formulas from $\rmvlog$  have variables from the set $\{v_1,\ldots,v_n\}$.
As in the previous section, we use the  function
\bce
$\sq :{\mathbb R}\to [0,1]$,  $\sq(x)=(x\vee 0)\wedge 1$ for any $x\in {\mathbb R}$.
\ece

\bfact  If $f_1$, $\ldots$, $f_k :[0,1]^n\to{\mathbb R}$  are continuous  piecewise linear functions and  $c_1$, $\ldots$, $c_k\in {\mathbb R}$ then

\noindent $\sum_{i=1}^k c_if_i$ is also a continuous piecewise linear function,  so $\sq\circ(\sum_{i=1}^k c_if_i)$ is in
$PL_n$.
By Theorem \ref{main}, there exists a formula
$\vp$ of $\rmvlog$ such that $\fct{\vp}=\sq\circ(\sum_{i=1}^k c_if_i)$. Therefore, we introduce the following definition.
\efact

\bdfn
Let $f_1$, $\ldots$, $f_k:[0,1]^n\to{\mathbb R}$  be continuous  piecewise linear functions. 
We say that $\vp$ is a {\em quasi-linear combination} of $f_1$, $\ldots$, $f_k$ whenever
\bcev $\fct{\vp}=\sq\circ(\sum_{i=1}^k c_if_i)$\ecev  for some $c_1, \ldots, c_k\in {\mathbb R}$.
We define
$\li{f_1, \ldots, f_k}$ as the subset of $Form(\rmvlog)$ that contains all the  quasi-linear combinations  of $f_1$, $\ldots$, $f_k$.    If
 $\vp_1$, $\ldots$, $\vp_k$ are formulas of $\rmvlog$ then  $\li{\fct{\vp_1}, \ldots, \fct{\vp_k}}$ will be denoted by
 $\li{\vp_1, \ldots, \vp_k}$.  If $\vp\in \li{\vp_1, \ldots, \vp_k}$ then we say that $\vp$ is a {\em quasi-linear combination of the formulas} $\vp_1$, $\ldots$, $\vp_k$ in $\rmvlog$.
\edfn

\blem\label{eval-lin}
If  $\vp$, $\vp_1$, $\ldots$, $\vp_k$ are formulas in $\rmvlog$ such that
$\fct{\vp}=\sq\circ(\sum_{i=1}^k c_i\fct{\vp_i})$ for some $c_1$, $\ldots$, $c_k\in {\mathbb R}$ and
$e:Form({\mathbb R}{\cal L})\sra [0,1]$ is an evaluation, then
\bce
$e(\vp)= \sq(\sum_{i=1}^k  c_i e(\vp_i))$.
\ece
\elem
\begin{proof} We have
$ e(\vp)=\fct{\vp}(e(v_1),\cdots, e(v_n))=$\\
$ \sq(\sum_{i=1}^k  c_i \fct{\vp_i}(e(v_1),\cdots, e(v_n)))=\sq(\sum_{i=1}^k  c_i e(\vp_i))$.  
\end{proof}

\blem
If  $\vp_1$, $\ldots$, $\vp_k$ are formulas of  ${\mathbb R}{\cal L}$  and\\   $r_1$, $\cdots$, $r_k\in [0,1]$ then
\bce
$\nabla_{r_1}\vp_1\oplus\cdots\oplus\nabla_{r_k}\vp_k \in \li{\vp_1, \ldots ,\vp_k}$.
\ece
\elem
\begin{proof}
Under the above hypothesis, we get
\bce
$\sq(\sum_{i=1}^k  r_i e(\vp_i))= r_1e(\vp_1)\oplus\cdots\oplus {r_k}e(\vp_k)=e(\nabla_{r_1}\vp_1)\oplus\cdots\oplus e(\nabla_{r_k}\vp_k)$
\ece
for any evaluation $e:Form({\mathbb R}{\cal L})\sra [0,1]$.  
\end{proof}

\bdfn
We say that a formula $\vp$ of ${\mathbb R}{\cal L}$ is {\em invalid} if there exists an evaluation
$e:Form({\mathbb R}{\cal L})\sra [0,1]$ such that $e(\vp)=0$
\edfn

\bthm\label{invth}
Let  $\vp_1$, $\ldots$, $\vp_k$ are formulas of  ${\mathbb R}{\cal L}$  and   $r_1$, $\cdots$, $r_k\in [0,1]$. The following are equivalent:\\
(i)  the book $\{(\vp_i,r_i)\mid i\in\{1,\ldots,k\}\}$ is {coherent},\\
(ii) there exists a state $s:RL_n\to [0,1]$ such that $s([\vp_i])=r_i$ for any $i\in\{1,\ldots,k\}$,\\
(iii) $\li{\fct{\vp_1}-r_1,\ldots, \fct{\vp_k}-r_k}$ is a set of invalid formulas of $\rmvlog$.
\ethm
\begin{proof}
(i)$\Leftrightarrow$(ii) Apply Corollary \ref{th:dfrmv} to $\beta([\vp_i])=r_i$ for any $i\in\{1,\ldots,k\}$.\\
(i)$\Leftrightarrow$(iii) The following facts are equivalent:\\
(1) $\li{\fct{\vp_1}-r_1,\ldots, \fct{\vp_k}-r_k}$ is a set of invalid formulas,\\
(2) for any $\Psi\in \li{\fct{\vp_1}-r_1,\ldots, \fct{\vp_k}-r_k}$ there exists an evaluation
$e:Form({\mathbb R}{\cal L})\sra [0,1]$ such that $e(\Psi)=0$,\\ 
(3) for any $c_1$, \ldots, $c_k\in {\mathbb R}$, if $\fct{\Psi}=\sq\circ (\sum_{i=1}^k c_i(\fct{\vp_i}-r_i))$ then there exists an evaluation
$e:Form({\mathbb R}{\cal L})\sra [0,1]$ such that $e(\Psi)=0$,\\
(4) for any $c_1$, \ldots, $c_k\in {\mathbb R}$, there exists an evaluation
$e:Form({\mathbb R}{\cal L})\sra [0,1]$ such that
\bce
$\sq\circ (\sum_{i=1}^k c_i(e(\vp_i)-r_i)=0$,
\ece
(5) for any $c_1$, \ldots, $c_k\in {\mathbb R}$, there exists an evaluation
$e:Form({\mathbb R}{\cal L})\sra [0,1]$ such that 
\bce $\sum_{i=1}^k c_i(r_i-e(\vp_i))\geq 0$.\ece
Note that  by Lemma \ref{eval-lin} is used for (3)$\Leftrightarrow$(4) and  Lemma \ref{calc} (b) is used 
 for (4)$\Leftrightarrow$(5). 
  
\end{proof}

\medskip

 If  $r\in [0,1]$ and $\vp$ a formula then the piecewise linear
  function $r-\fct{\vp}: [0,1]^n\to{\mathbb R}$  may have  negative values, therefore  it may not correspond to a formula of ${\mathbb R}{\cal L}$.
 The next result provides
 a necessary condition for a book to be coherent using quasi-linear combinations of formulas.  We also prove a sufficient condition in  Corollary \ref{invpr2}, but using different formulas.

\medskip

 We set  $\bm{r}=\Delta_r(\varphi\sra\varphi)$  and $\vp\ominus\psi=\vp\odot\neg\psi$ whenever $r\in [0,1]$ and $\varphi$, $\psi\in Form({\mathbb R}{\cal L})$. Note that $\chi=\vp\ominus \psi$ implies $\fct{\chi}=0\vee(\fct{\vp}-\fct{\psi})$.

 \medskip

\bprop\label{invpr1}

Assume   $\vp_1$, $\ldots$, $\vp_k$ are formulas of  ${\mathbb R}{\cal L}$  and 
$r_1$, $\cdots$, $r_k\in [0,1]$ such that
the book $\{(\vp_i,r_i)\mid i\in\{1,\ldots,k\}\}$ is {coherent}. Hence
\bce
$\li{(\bm{r_1}\ominus\vp_1),\ldots, (\bm{r_k}\ominus\vp_k) }$ 
\ece
is  a set of invalid formulas.
\eprop
\begin{proof} Let $\chi\in \li{(\bm{r_1}\ominus\vp_1),\ldots, (\bm{r_k}\ominus\vp_k) }$ and\\
 $c_1$, $\cdots$, $c_k\in {\mathbb R}$  such that
\bce
$\fct{\chi}=\sq\circ(\sum_{i=1}^kc_i(0\vee(r_i-\fct{\vp_i})))$.
\ece

Since the book  $\{(\vp_i,r_i)\mid i\in\{1,\ldots,k\}\}$ is {coherent} there is an evaluation $e$ such that
\bce
$\sum_{i=1}^k(-c_i)(r_i-e(\vp_i))\geq 0$, which implies that $\sum_{i=1}^k(-c_i)((r_i-e(\vp_i))\vee 0)\geq 0$.
\ece

We get $\sum_{i=1}^k(-c_i) e(\bm{r_i}\ominus\vp_i)\geq 0$.  By Lemma \ref{calc} (b), $\sq(\sum_{i=1}^k c_i e(\bm{r_i}\ominus\vp_i))=0$.
Using  Lemma \ref{eval-lin} it follows that  $e(\chi)=0$, so $\chi$ is an invalid formula.  
\end{proof}

If $\vp\in Form({\mathbb R}{\cal L})$, $r\in [0,1]$ and $c\in {\mathbb R}$ we denote
\bce
$\psi(\vp,r,c)=\left\{\begin{array}{ll}
\vp\ominus \bm{r}, & \mbox{ if } c\geq 0\\
\bm{r}\ominus\vp, &  \mbox{ if } c < 0.
\end{array}\right.$
\ece

\bprop\label{prajut}
Assume   $\vp_1,\ldots,\vp_k\in Form({\mathbb R}{\cal L})$  and   $r_1$, $\cdots$, $r_k\in [0,1]$ such that
for any  $c_1,\ldots,c_k\in {\mathbb R}$   the formula $\Phi$  of $\rmvlog$ is invalid whenever
$\fct{\Phi}= \sq\circ(\sum_{i=1}^k|c_i|\fct{\psi_i})$,
 where $|c_i|$ is the module of $c_i$ and $\psi_i=\psi(\vp_i,r_i,c_i)$ for any $i\in\{1,\ldots,k\}$.
Then the book $\{(\vp_i,r_i)\mid i\in\{1,\ldots,k\}\}$ is {coherent}.
\eprop
\begin{proof} If $c_1$, $\cdots$, $c_k\in {\mathbb R}$ and   $\Phi$ in $Form(\rmvlog)$  is an invalid formula 
such that  $\fct{\Phi}= \sq\circ(\sum_{i=1}^k|c_i|\fct{\psi_i})$,
 then there exists an evaluation  $e$ such that $e(\Phi)=0$ and, by \mbox{ Lemma \ref{eval-lin},}  we get
$\sq(\sum_{i=1}^k  |c_i| e(\psi_i))=0$, so
\bce
$\begin{array}{ccc}
\sq\left( \right. & \sum_{c_i\geq 0}  c_i((e(\vp_i)-r_i)\vee 0) +& \\

 & \sum_{c_i< 0} (-c_i)((r_i-e(\vp_i))\vee 0) & \left. \right)=0
\end{array}$
\ece
By Lemma \ref{calc}(a), 
\bce
 $\sq((\sum_{c_i\geq 0} c_i(e(\vp_i)-r_i) +\sum_{c_i< 0} (-c_i)(r_i-e(\vp_i)))\vee 0)=0$
\ece
 and,  by Lemma \ref{calc}(c),
\bce
$\sq(\sum_{c_i\geq 0} c_i(e(\vp_i)-r_i) +\sum_{c_i< 0} (-c_i)(r_i-e(\vp_i)))=0$.
\ece
Using  Lemma \ref{calc} (b), we get
\bce
$-(\sum_{c_i\geq 0} c_i(e(\vp_i)-r_i) +\sum_{c_i< 0} (-c_i)(r_i-e(\vp_i)))\geq 0$.
\ece
It follows that  $\sum_{i=1}^k  c_i (r_i-e(\vp_i))\geq 0$.    
\end{proof}

\medskip

\bcor\label{invpr2}
Assume   $\vp_1$, $\ldots$, $\vp_k\in Form({\mathbb R}{\cal L})$  and \\  $r_1$, $\cdots$, $r_k\in [0,1]$ such that
\bce
$\li{\alpha_1,\ldots,\alpha_k,\beta_1,\ldots,\beta_k}$
\ece
 is a set of invalid formulas, where 
\bce
$\alpha_i=\bm{r_i}\ominus\vp_i$ and $\beta_i=\vp_i\ominus \bm{r_i}$ for any  $i\in\{1,\ldots,k\}$.
\ece 
Then the book $\{(\vp_i,r_i)\mid i\in\{1,\ldots,k\}\}$ is {coherent}.
\ecor
\begin{proof}
For  any  $c_1$, $\cdots$, $c_k\in {\mathbb R}$, if $\Phi$ is a formula of $\rmvlog$ such that $\fct{\Phi}= \sq\circ\sum_{i=1}^k|c_i|\psi_i$ as in
Proposition \ref{prajut}, then $\Phi\in \li{\alpha_1,\ldots,\alpha_k,\beta_1,\ldots,\beta_k}$. In consequence, we can apply Proposition \ref{prajut}.  
\end{proof}

\bfact
We initiate the theory of  quasi-linear combinations in ${\mathbb R}{\cal L}$ and we relate it with  de Finetti's notion of coherence, which can be expressed by an invalidity condition.
 The linear combinations of formulas from \L ukasiewicz logic were approached in \cite{Trillas} and \cite{Amato}, as representations for a particular class of neural networks.
The composition between the  function $\sq$ and a linear combination of formulas \L ukasiewicz logic  can be naturally represented by a formula in
our logic ${\mathbb R}{\cal L}$ and, therefore, the theory of linear combinations can be approached within a simple defined logical system. Note that ${\mathbb R}{\cal L}$ is a conservative extension of \L ukasiewicz  logic.
It has standard completeness theorem with respect to $[0,1]$ and it is supported by the algebraic theory of Riesz MV-algebras which are categorically equivalent with unital Riesz spaces.  Hence, our hope for the future is that
the system  ${\mathbb R}{\cal L}$ is enough expressive for representing classes of neural networks in a pure logical frame, having in mind the role of classical logic in the synthesis and analysis of boolean circuits.
\efact

\bigskip

\noindent{\bf Acknowledgement.}\\
I. Leu\c stean  was supported by the strategic grant POSDRU/89/1.5/S/58852,  cofinanced by ESF within SOP HRD 2007-2013. Part of the research has been carried out while visiting the University of Salerno. 
The investigations from Section \ref{sec-norm} were initiated in 2009, when V. Marra (University of Milan) pointed out Kakutani's theorem and its  relevance to the theory of Riesz MV-algebras.

\end{document}